\documentclass{article}
 \usepackage{amssymb,amsmath,amsfonts}
 -\marginparwidth \oddsidemargin -4mm \evensidemargin
 \oddsidemargin%
 \advance\hoffset by 5mm

 \title{Asymptotic stability of ground states in 2D nonlinear Schr\" odinger equation including subcritical cases}
 \author{E. Kirr \thanks{Department of Mathematics, University of Illinois at Urbana-Champaign}
 \hspace{0.05in} and A. Zarnescu \thanks{Department of Mathematics, University of Oxford, England}}
 \date{\today}
 \newtheorem{lemma}{Lemma}[section]
 \newtheorem{theorem}{Theorem}[section]
 \newtheorem{proposition}{Proposition}[section]
 \newtheorem{corollary}{Corollary}[section]
 \newtheorem{remark}{Remark}[section]

\begin{document}

 \maketitle
 \begin{abstract}
 \noindent We consider a class of nonlinear Schr\"{o}dinger equation in two space
 dimensions with an attractive potential. The nonlinearity is local but rather general
 encompassing for the first time both subcritical and supercritical (in $L^2$) nonlinearities. We study the asymptotic
 stability of the nonlinear bound states, i.e. periodic in time
 localized in space solutions. Our result shows that all solutions
 with small initial data, converge to a nonlinear bound state.
 Therefore, the nonlinear bound states are asymptotically stable. The proof hinges on dispersive estimates that we
 obtain for the time dependent, Hamiltonian, linearized dynamics
 around a careful chosen one parameter family of bound states that ``shadows" the nonlinear evolution of the
 system. Due to the generality of the methods we develop we expect
 them to extend to the case of perturbations of large bound states and to other nonlinear dispersive wave type equations.
 \end{abstract}

\makeatletter \@addtoreset{equation}{section} \makeatother
\renewcommand{\theequation}{\thesection.\arabic{equation}}

\section{Introduction}

In this paper we study the long time behavior of solutions of the
nonlinear Schr\" odinger equation (NLS) with potential in two space
dimensions (2-d):
\begin{eqnarray}
i\partial_t u(t,x)&=&[-\Delta_x+V(x)]u+g(u), \quad t\in\mathbb{R},\quad x\in\mathbb{R}^2\label{eq:ufull}\\
u(0,x)&=&u_0(x)\label{eq:ic}
\end{eqnarray}
where the local nonlinearity is constructed from the real valued,
odd, $C^2$ function $g:\mathbb{R}\mapsto\mathbb{R}$ satisfying
\begin{equation}\label{gest}
g(0)=g'(0)=0\quad {\rm and}\quad |g''(s)|\leq
C(|s|^{\alpha_1}+|s|^{\alpha_2}),\quad s\in\mathbb{R},\
\frac{1}{2}<\alpha_1\le\alpha_2<\infty\end{equation} which is then
extended to a complex function via the gauge symmetry:
\begin{equation}\label{gsym}
g(e^{i\theta}s)=e^{i\theta}g(s),\qquad \theta\in\mathbb{R}.
\end{equation}
The equation has important applications in statistical physics,
optics and water waves. It describes certain limiting behavior of
Bose-Einstein condensates \cite{dgps:bec,lsy:2d} and propagation of
time harmonic waves in wave guides \cite{kn:kp,kn:Marcuse,nm:no}. In
the latter, $t$ plays the role of the coordinate along the axis of
symmetry of the wave guide.

It is well known that this nonlinear equation admits periodic in
time, localized in space solutions (bound states or solitary waves).
They can be obtained via both variational techniques
\cite{bl:i,str:sw,rw:bs} and bifurcation methods
\cite{pw:cm,rw:bs,kz:as2d}. Moreover the set of periodic solutions
can be organized as a $C^2$ manifold (center manifold), see
\cite{gnt:as,km:as3d} or next section. Orbital stability of solitary
waves, i.e. stability modulo the group of symmetries $u\mapsto
e^{-i\theta}u,$ was first proved in \cite{rw:bs,mw:ls}, see also
\cite{gss:i,gss:ii,ss:ins}.

The main result of this paper is that solutions of
\eqref{eq:ufull}-\eqref{eq:ic} with small initial data
asymptotically converge to a bound state, see Theorem \ref{th:main}.
While asymptotic stability results for bound states in NLS have
first appeared in the work of A. Soffer and M. I. Weinstein
\cite{sw:mc1,sw:mc2}, and continued in
\cite{pw:cm,kn:Wed,bp:asi,bp:asii,bs:as,sc:as,gnt:as}, our main
contribution is to allow for subcritical and critical ($L^2$)
nonlinearities, $\frac{1}{2}<\alpha_1\le 1$ in \eqref{gest}. To
accomplish this we carefully project the nonlinear dynamics onto the
center manifold of bound states and use linearization around this
time changing projection to study the motion in the radiative
directions, i.e. directions that are not in the tangent space of the
center manifold. Previously, linearization around a fixed bound
state has been used, see the papers cited above. By continuously
adapting the linear dynamics to the actual nonlinear evolution of
the solution we can more precisely capture the effective potential
induced by the nonlinearity $g$ into a time dependent linear
operator. Once we have a good understanding of this time dependent
linear dynamics, i.e. we have good dispersive estimates for its
semigroup of operators, see Section \ref{se:lin}, we obtain
information for the nonlinear dynamics via Duhamel formula and
contraction principles for integral equations, see Section
\ref{se:main}. Note that we have recently used a similar technique
to show that in the critical (cubic) case, \eqref{eq:ufull} with
$g(s)=s^3,\ s\in\mathbb{R},$ the center manifold of bound states is
an attractor for small initial data, see \cite{kz:as2d}. In this
paper the technique is much refined, we use a better projection of
the dynamics on the center manifold and sharper estimates for the
linear dynamics. The refinements not only allow us to treat a much
larger spectrum of nonlinearities including, for the first time, the
subcritical ones but also allow us to obtain actual convergence of
the solution to a bound state.

However, the main challenge for our approach is to obtain good
dispersive estimates for the semigroup of operators generated by the
time dependent linearization that we use. This is accomplished in
Section \ref{se:lin} via a perturbative method similar to the one we
developed in \cite{kz:as2d}. As described in that section, we could
have obtained sharper estimates by using a generalized Fourier
multiplier technique to remove the singularity of
$$\|e^{i(\Delta -V) t}\|_{L^{1}\mapsto L^\infty}\sim |t|^{-1},$$
see \cite[Section 4]{km:as3d}. We chose not to do it because it
requires stronger hypotheses on $V$ without allowing us to enlarge
the spectrum of nonlinearities that we can treat.

Finally, we remark that our method is quite general, based solely on
linearization around nonlinear bound states and estimates for
integral operators with dispersive kernels. Therefore we expect it
to generalize to the case of large nonlinear 2D ground states, see
for example \cite{sc:as}, the presence of multiple families of bound
states, see for example \cite{sw:sgs}, or to the case of time
dependent nonlinearity, see \cite{ckp:res}. In all three cases our
method will not only allow to treat the less dispersive environment,
2D compared to 3D, but it should greatly reduce the restrictions on
the nonlinearity. The first author and collaborators are currently
working on adapting the method to other dimensions and other
dispersive wave type equations. The work in 3-D is complete, see
\cite{km:as3d}.

\bigskip

\noindent{\bf Notations:} $H=-\Delta+V;$

$L^p=\{f:\mathbb{R}^2\mapsto \mathbb{C}\  |\ f\ {\rm measurable\
and}\ \int_{\mathbb{R}^2}|f(x)|^pdx<\infty\},\ 1\le p<\infty,$
endowed with the standard norm
$\|f\|_{L^p}=\left(\int_{\mathbb{R}^2}|f(x)|^pdx\right)^{1/p},$
while for $p=\infty,$ $L^\infty=\{f:\mathbb{R}^2\mapsto \mathbb{C}\
|\ f\ {\rm measurable\ and}\ {\rm essup}|f(x)|<\infty\},$ and it is
endowed with the norm: $\|f\|_{L^\infty}={\rm essup}|f(x)|$;

$<x>=(1+|x|^2)^{1/2},$ and for $\sigma\in\mathbb{R},\ 1\le
p<\infty,$ $L^p_\sigma$ denotes the $L^p$ space with weight
$<x>^{p\sigma},$ i.e. the space of functions $f(x)$ such that
$(<x>^{\sigma}f(x))^p$ are integrable endowed with the norm
$\|f(x)\|_{L^p_\sigma}=\|<x>^{\sigma}f(x)\|_p,$ while for
$p=\infty,$ $L^\infty_\sigma$ denotes the vector space of measurable
functions $f(x)$ such that ${\rm essup}|<x>^{\sigma}f(x)|<\infty$
endowed with the norm
$\|f(x)\|_{L^\infty_\sigma}=\|<x>^{\sigma}f(x)\|_{L^\infty};$

$\langle f,g\rangle =\int_{\mathbb{R}^2}\overline f(x)g(x)dx$ is the
scalar product in $L^2$ where $\overline z=$ the complex conjugate
of the complex number $f;$

$P_c$ is the projection associated to the continuous spectrum of the
self adjoint operator $H$ on $L^2,\ {\rm range} P_c={\cal H}_0;$

$H^n$ denote the Sobolev spaces of measurable functions having all
distributional partial derivatives up to order $n$ in $L^2,
\|\cdot\|_{H^n}$ denotes the standard norm in this spaces.

\section{Preliminaries. The center manifold.}\label{se:prelim}

The center manifold is formed by the collection of periodic solutions for
(\ref{eq:ufull}):
\begin{equation}\label{eq:per}
  u_E(t,x)=e^{-iEt}\psi_E(x)
\end{equation}
where $E\in\mathbb{R}$ and $0\not\equiv\psi_E\in H^2(\mathbb{R}^2)$
satisfy the time independent equation:
\begin{equation}\label{eq:ev}
[-\Delta+V]\psi_E+g(\psi_E)=E\psi_E
\end{equation}
Clearly the function constantly equal to zero is a solution of
(\ref{eq:ev}) but (iii) in the following hypotheses on the potential
$V$ allows for a bifurcation with a nontrivial, one (complex)
parameter family of solutions:

\bigskip\noindent{\bf (H1)} Assume that
\begin{itemize}
  \item[(i)] There exists $C>0$ and $\rho >3$ such that:
  $$|V(x)|\le C<x>^{-\rho},\quad {\rm for\ all}\ x\in\mathbb{R}^2;$$
  \item[(ii)] $0$ is a regular point\footnote{see
 \cite[Definition 7]{ws:de2} or $M_\mu=\{0\}$ in relation (3.1) in \cite{mm:ae}}
 of the
spectrum of
  the linear operator $H=-\Delta+V$ acting on $L^2;$
  \item [(iii)]$H$ acting on $L^2$ has exactly one
  negative eigenvalue $E_0<0$ with corresponding normalized
  eigenvector $\psi_0.$ It is well known that $\psi_0(x)$ can be
  chosen strictly positive and exponentially decaying as
$|x|\rightarrow\infty.$
\end{itemize}

\par\noindent Conditions (i)-(ii) guarantee the applicability of dispersive
estimates of Murata \cite{mm:ae} and Schlag \cite{ws:de2} to the
Schr\" odinger group $e^{-iHt}.$ These estimates are used for
obtaining Theorems \ref{th:lin1} and \ref{th:lin2}, see also
\cite[section 4]{kz:as2d}. In particular (i) implies the local well
posedness in $H^1$ of the initial value problem
(\ref{eq:ufull}-\ref{eq:ic}), see section \ref{se:main}.

By the standard bifurcation argument in Banach spaces \cite{ln:fa}
for (\ref{eq:ev}) at $E=E_0,$ condition (iii) guarantees existence
of nontrivial solutions. Moreover, these solutions can be organized
as a $C^2$ manifold (center manifold) for $x\in\mathbb{R}^n,$ see
\cite[section 2]{km:as3d} or \cite{gnt:as}. The proofs for the
following results can be found in \cite[section 2]{km:as3d} or
\cite{gnt:as}:

\begin{proposition}\label{pr:cm} There exist $\delta>0,$
the $C^2$ function
$$h:\{a\in\mathbb{R}\times\mathbb{R}\ :\ |a|<\delta\}\mapsto L^2_\sigma\cap H^2
,\ \sigma\in\mathbb{R}$$ and the $C^1$ function
$E:(-\delta,\delta)\mapsto\mathbb{R}$ such that for $|E-E_0|<\delta$
and $|\langle\psi_0,\psi_E\rangle |<\delta$ the eigenvalue problem
(\ref{eq:ev}) has a unique solution up to multiplication with
$e^{i\theta},\ \theta\in [0,2\pi),$ which can be represented as a
center manifold:
 \begin{equation}\label{eq:cm}
 \psi_E=a\psi_0+h(a),\ E=E(|a|), \quad \langle\psi_0,h(a)\rangle =0,\quad
 h(e^{i\theta}a)=e^{i\theta}h(a),\
|a|<\delta .\end{equation} Moreover
$E(|a|)=\mathcal{O}(|a|^{1+\alpha_1})$,
$h(a)=\mathcal{O}(|a|^{2+\alpha_1}),$ and for $a\in\mathbb{R},\
|a|<\delta,$ $h(a)$ is a real valued function with
$\frac{d^2h}{da^2}(a)=\mathcal{O}(|a|^{\alpha_1})$ and
$\frac{dh}{da}(0)=0.$
\end{proposition}

Since $\psi_0(x)$ is exponentially decaying as
$|x|\rightarrow\infty$ the proposition implies that $\psi_E\in
L^2_\sigma .$ A regularity argument, see \cite{sw:mc1}, gives a
stronger result:

\begin{corollary}\label{co:decay} For any $\sigma\in\mathbb{R},$
there exists a finite constant $C_\sigma$ such that:
$$\|<x>^\sigma\psi_E\|_{H^2}\le C_\sigma\|\psi_E\|_{H^2}.$$
\end{corollary}

We are going to decompose the solution of
\eqref{eq:ufull}-\eqref{eq:ic} into a projection onto the center
manifold and a correction. To insure that the correction disperses
to infinity on long times we require that the correction is always
in the invariant subspace of the linearized dynamics at the
projection that complements the tangent space to the center
manifold. A short description of the decomposition follows, for more
details and the proofs see \cite{km:as3d}.

Consider the linearization of \eqref{eq:ufull} at a function on the
center manifold $\psi_E=a\psi_0+h(a),\ a=a_1+ia_2\in\mathbb{C},\
|a|<\delta:$
 \begin{equation}\label{eq:ldE}
 \frac{\partial w}{\partial t}=-iL_{\psi_E}[w]-iEw
 \end{equation}
where
 \begin{equation}\label{def:linop}
 L_{\psi_E}[w]=(-\Delta+V-E)w+Dg_{\psi_E}[w]=(-\Delta+V-E)w+\lim_{\varepsilon\in\mathbb{R},\ \varepsilon\rightarrow
 0}\frac{g(\psi_E+\varepsilon w)-g(\psi_E)}{\varepsilon}
 \end{equation}

\begin{remark}\label{rmk:dgest} Note that for $a\in\mathbb{R}$ we
have $\psi_E=a\psi_0+h(a)$ is real valued and
 $$Dg_{\psi_E}[w]=g'(\psi_E)\Re w+i\frac{g(\psi_E)}{\psi_E}\Im w
   =\frac{1}{2}\left(g'(\psi_E)+\frac{g(\psi_E)}{\psi_E}\right)w
   +\frac{1}{2}\left(g'(\psi_E)-\frac{g(\psi_E)}{\psi_E}\right)\overline{w}$$
hence
 \begin{equation}\label{est:dg}
 |Dg_{\psi_E}[w]|\le |w|\max\left\{|g'(\psi_E)|,
 \left|\frac{g(\psi_E)}{\psi_E}\right|\right\}\le
 C(|\psi_E|^{1+\alpha_1}+|\psi_E|^{1+\alpha_2})|w|
 \end{equation}
where we used \eqref{gest}. For $a=|a|e^{i\theta}\in\mathbb{C}$ we
have, using the equivariant symmetry \eqref{gsym},
$\psi_E=a\psi_0+h(a)=e^{i\theta}(|a|\psi_0+h(|a|)=e^{i\theta}\psi_E^{\rm
real},$ where $\psi_E^{\rm real}$ is real valued, and
$Dg_{\psi_E}[w]=e^{i\theta}Dg_{\psi_E^{\rm real}}[e^{-i\theta}w],$
hence \eqref{est:dg} is valid for any $\psi_E$ on the manifold of
ground states. \end{remark}

{\bf Properties of the linearized operator}:
 \begin{enumerate}
 \item $L_{\psi_E}$ is real linear and symmetric with respect to the
 real scalar product $\Re\langle\cdot,\cdot\rangle,$ on
 $L^2(\mathbb{R}^2),$ with domain $H^2(\mathbb{R}^2).$
 \item Zero is an e-value for $-iL_{\psi_E}$ and its generalized
 eigenspace includes $\left\{\frac{\partial\psi_E}{\partial a_1},\frac{\partial\psi_E}{\partial
 a_2}\right\}$
 \item ${\rm span}_{\mathbb{R}}\left\{\frac{\partial\psi_E}{\partial a_1},\frac{\partial\psi_E}{\partial
 a_2}\right\}$ and
  ${\cal H}_a=\left\{-i\frac{\partial\psi_E}{\partial a_2},i\frac{\partial\psi_E}{\partial
 a_1}\right\}^\perp,$
 where orthogonality is with respect to the real scalar product in $L^2(\mathbb{R}^2),$ are
 invariant subspaces for $-iL_{\psi_E}$ and, by possible choosing
 $\delta>0$ smaller than the one in Proposition \ref{pr:cm}, we
 have:
 $$L^2(\mathbb{R}^2)={\rm span}_{\mathbb{R}}\left\{\frac{\partial\psi_E}{\partial a_1},\frac{\partial\psi_E}{\partial
 a_2}\right\}\oplus {\cal H}_a,\qquad {\rm for\ all}\ a\in\mathbb{C},\ |a|<\delta.$$ Note that ${\cal H}_0$ coincides
 with the subspace of $L^2$ associated to the continuous spectrum of
 the self-adjoint operator
 $H=-\Delta+V.$

 \item the above decomposition can be extended to
 $H^{-1}(\mathbb{R}^2):$
 \begin{equation}\label{h-1decomp}H^{-1}(\mathbb{R}^2)={\rm span}_{\mathbb{R}}\left\{\frac{\partial\psi_E}{\partial a_1},\frac{\partial\psi_E}{\partial
 a_2}\right\}\oplus {\cal H}_a,\qquad {\rm for\ all}\ a\in\mathbb{C},\ |a|<\delta,\end{equation}
 where
 $${\cal H}_a=\left\{\phi\in H^{-1}\ |\ \Re\langle -i\frac{\partial\psi_E}{\partial
 a_2},\ \phi\rangle=0,\ {\rm and}\ \Re\langle i\frac{\partial\psi_E}{\partial
 a_1},\ \phi\rangle=0\right\}$$

 \end{enumerate}

Our goal is to decompose the solution of \eqref{eq:ufull} at each
time into:
 $$u=\psi_E+\eta=a\psi_0+h(a)+\eta,\qquad \eta\in{\cal H}_a$$
which insures that $\eta$ is not in the non-decaying directions of
the linearized equation \eqref{eq:ldE} at $\psi_E.$ The fact that
this can be done in an unique manner is a consequence of the
following lemma:

\begin{lemma}\label{lem:decomp} There exists $\delta /2>\delta_1>0$ such that
 any $\phi\in H^{-1}(\mathbb{R}^2)$ satisfying $\|\phi\|_{H^{-1}}\le\delta_1$
can be uniquely decomposed:
 $$\phi =\psi_E+\eta=a\psi_0+h(a)+\eta$$
where $a=a_1+ia_2\in\mathbb{C},\ |a|<\delta,\ \eta\in {\cal H}_a.$
Moreover the maps $\phi\mapsto a$ and $\phi\mapsto \eta$ are $C^1$
and there exists the constant $C$ independent on $\phi$  such that
$$|a|\le 2\|\phi\|_{H^{-1}},\qquad \|\eta\|_{H^{-1}}\le C\|\phi\|_{H^{-1}},$$
while for $\phi\in L^2(\mathbb{R}^2)$ we have $\eta\in
L^2(\mathbb{R}^2)$ and:
$$|a|\le 2\|\phi\|_{L^2},\qquad \|\eta\|_{L^2}\le C\|\phi\|_{L^2}.$$
\end{lemma}

\begin{remark}\label{rmk:inv} The above lemma uses the implicit
function theorem applied to $$F:\mathbb{R}^2\times
H^{-1}(\mathbb{R}^2)\mapsto\mathbb{R}^2\qquad
 F(a_1,a_2,\phi)=\left[\begin{array}{c}\Re\langle \Psi_1,\ \psi_E-\phi\rangle\\
 \Re\langle \Psi_2,\ \psi_E-\phi\rangle\end{array}\right]$$
where $\psi_E=(a_1+ia_2)\psi_0+h(a_1+ia_2)$ and
 \begin{eqnarray}
 \Psi_1(a_1,a_2)&=&-i\frac{\partial\psi_E}{\partial
 a_2}\left(\Re\langle -i\frac{\partial\psi_E}{\partial
 a_2},\ \frac{\partial\psi_E}{\partial
 a_1}\rangle\right)^{-1}\nonumber\\
 \Psi_2(a_1,a_2)&=&i\frac{\partial\psi_E}{\partial
 a_1}\left(\Re\langle i\frac{\partial\psi_E}{\partial
 a_1},\ \frac{\partial\psi_E}{\partial
 a_2}\rangle\right)^{-1}\nonumber
 \end{eqnarray}
form the dual basis of $\left\{\frac{\partial\psi_E}{\partial
a_1},\frac{\partial\psi_E}{\partial a_2}\right\}$ with respect to
the decomposition \eqref{h-1decomp}. Note that
$$\frac{\partial F}{\partial
(a_1,a_2)}(a_1,a_2,\phi)=\mathbb{I}_{\mathbb{R}^2}-M(a_1,a_2,\phi)$$
where the entries of the two by two matrix $M$ are
$$M_{ij}=\Re\langle \frac{\partial\Psi_i}{\partial a_j},\
\phi-\psi_E\rangle$$ and, consequently, $M(0,0,0)$ is the zero
matrix. Thus the implicit function theorem applies to $F=0,$ in a
neighborhood of $(a_1,a_2,\phi)=(0,0,0)$ and the number $\delta_1$
in the above lemma is chosen such that:
 $$\left|\Re\langle i\frac{\partial\psi_E}{\partial
 a_1},\ \frac{\partial\psi_E}{\partial
 a_2}\rangle\right|\ge \frac{1}{2},\qquad {\rm whenever}\
 |(a_1,a_2)|\le 2\delta_1,$$
and the norm of the matrix $M$ as a linear, bounded operator on
$\mathbb{R}^2$ satisfies:
\begin{equation}\label{M-bound}
\|M_\phi\|=\|M(a_1(\phi),a_2(\phi),\phi)\|\le\frac{1}{2},\qquad {\rm
whenever}\ \|\phi\|_{H^{-1}}\le\delta_1,
\end{equation}
see \cite[section 2]{km:as3d} for details.
\end{remark}

We need one more technical result relating the spaces ${\cal H}_a$
and the space corresponding to the continuous spectrum of
$-\Delta+V:$
 \begin{lemma}\label{le:pcinv}
With $\delta_1$ given by the previous lemma we have that for any
$a\in\mathbb{C},\ |a|\le 2\delta_1,$ the linear map $P_c|_{{\cal
H}_a}:{\cal H}_a\mapsto {\cal
 H}_0$ is invertible, and its inverse $R_a :{\cal H}_0\mapsto {\cal
 H}_a$ satisfies:
 \begin{eqnarray}
 \|R_a\zeta\|_{L^2_{-\sigma}}&\le
 &C_{-\sigma}\|\zeta\|_{L^2_{-\sigma}},\qquad \sigma\in\mathbb{R}\ {\rm and\ for\ all}\ \zeta\in {\cal H}_0\cap L^2_{-\sigma}\nonumber\\
 \|R_a\zeta\|_{L^p}&\le
 &C_p\|\zeta\|_{L^p},\qquad 1\le p\le\infty\ {\rm and\ for\ all}\ \zeta\in {\cal H}_0\cap
 L^p\nonumber\\
 \overline{R_a\zeta}&=&R_a\overline\zeta\nonumber
 \end{eqnarray}
where the constants $C_{-\sigma},\ C_p>0$ are independent of
$a\in\mathbb{C},\ |a|\le 2\delta_1.$
 \end{lemma}

We are now ready to prove our main result.

\section{The Main Result}\label{se:main}

 \begin{theorem}\label{th:main} If hypothesis \eqref{gest}, \eqref{gsym}, $(H1)$ hold and
$$\frac{1}{2}<\alpha_1$$
then there
 exists $q'_0<\frac{4+2\alpha_2}{3+2\alpha_2}$ and $\varepsilon_0>0$ such that for all initial conditions
 $u_0(x)$ satisfying
 $$\max\{\|u_0\|_{L^{q'_0}},\|u_0\|_{H^1}\}\le \varepsilon_0
 $$ the initial value problem
(\ref{eq:ufull})-(\ref{eq:ic}) is
 globally well-posed in $H^1,$ and the solution decomposes into a radiative
part and a part that asymptotically converges to a ground state.

More precisely, there exist a $C^1$ function
$a:\mathbb{R}\mapsto\mathbb{C}$ such that, for all $t\in\mathbb{R}$
we have:
$$
u(t,x)=\underbrace{a(t)\psi_0(x)+h(a(t))}_{\psi_E(t)}+\eta(t,x) $$
where $\psi_E(t)$ is on the central manifold (i.e it is a ground
state) and $\eta(t,x)\in {\cal H}_{a(t)},$ see Proposition
\ref{pr:cm} and Lemma \ref{lem:decomp}. Moreover, there exists the
ground states $\psi_{E_{\pm\infty}}$ and the $C^1$ function
$\tilde\theta:\mathbb{R}\mapsto\mathbb{R}$ such that
$\lim_{|t|\rightarrow\infty}\theta(t)=0$ and:
 \begin{equation}\label{conv:psie}\lim_{t\rightarrow\pm\infty}\|\psi_E(t)-e^{-it(E_\pm-\theta(t))}\psi_{E_{\pm\infty}}\|_{H^2\bigcap
 L^2_\sigma}=0,\ \sigma\in\mathbb{R}
 \end{equation} while $\eta$ satisfies the
following decay estimates. Fix $p_0>\max\{\frac{2}{\alpha_1-1/2},\
(4+2\alpha_2)\frac{q_0-2}{q_0-(4+2\alpha_2)}\},$ where
$q_0=\frac{q'_0}{q'_0-1}>4+2\alpha_2.$ Then for $2\le
p\le\frac{p_0q_0}{p_0+q_0-2}$ we have:
 \begin{equation}\label{lpdecay}
 \|\eta(t)\|_{L^p}\le \left\{\begin{array}{cl}
     C\varepsilon_0\frac{\log^{\frac{1-2/p}{1-2/p_0}}(2+|t|)}{(1+|t|)^{1-2/p}}
     & {\rm if}\ \alpha_1\ge 1\ {\rm or}\ \alpha_1<1\ {\rm and}\ \ p\le \frac{2}{1-\alpha_1+2/p_0}, \\
      & \\
     C\varepsilon_0\frac{\log^{\frac{\alpha_1-2/p_0}{1-2/p_0}}(2+|t|)}{(1+|t|)^{\alpha_1-2/p_0}}
     & {\rm if}\ \alpha_1<1\ {\rm and}\ p>
     \frac{2}{1-\alpha_1+2/p_0},
     \end{array}\right.
 \end{equation}
for some constant $C=C(p_0).$
 \end{theorem}
\begin{remark}\label{rmk:rad} The estimates on $\eta$ show
that the component of the solution that does not converge to a
ground states disperses like the solution of the free Schr\" odinger
equation except for a logarithmic correction in $L^p$ spaces for
critical and supercritical regimes, $\alpha_1\ge 1.$ In subcritical
regimes, $\alpha_1<1,$ the decay rate remains comparable to the free
Schr\" odinger one in $L^p$ spaces for $2\le p <2/(1-\alpha_1),$
while it saturates to $|t|^{\alpha_1-1-0}$ in $L^p,\ p\ge
2/(1-\alpha_1).$
\end{remark}

\smallskip\par{\bf Proof of Theorem \ref{th:main}.} It is well
known that under hypothesis $(H1)(i)$ the initial value problem
(\ref{eq:ufull})-(\ref{eq:ic}) is locally well posed in the energy
space $H^1$ and its $L^2$ norm is conserved, see for example \cite[
Corollary 4.3.3. at p. 92]{caz:bk}. Global well posedness follows
via energy estimates from $\|u_0\|_{H^1}$ small, see \cite[Corollary
6.1.5 at p. 165]{caz:bk}.

We choose $\varepsilon_0\le \delta_1$ given by Lemma
\ref{lem:decomp}. Then, for all times, $\|u(t)\|_{H^{-1}}\le
\|u(t)\|_{L^2}\le\varepsilon_0\le\delta_1$ and, via Lemma
\ref{lem:decomp}, we can decompose the solution into a solitary wave
and a dispersive component:
\begin{equation}\label{dc}
u(t)=a(t)\psi_0+h(a(t))+\eta(t)=\psi_E(t)+\eta(t),\qquad {\rm
where}\ |a(t)|=|a_1(t)+ia_2(t)|\le 2\varepsilon_0\le 2\delta_1\
\forall t\in\mathbb{R}.
\end{equation} Note that since $a\mapsto h(a)$ is $C^2,$ see Proposition \ref{pr:cm}, and $a$ is uniformly bounded in
time we deduce that there exists the constant $C_H>0$ such that:
$$\max\left\{\|\psi_E(t)\|_{H^2},\|\frac{\partial\psi_E}{\partial
a_1}(t)\|_{H^2},\|\frac{\partial\psi_E}{\partial
a_2}(t)\|_{H^2}\right\}\le C_H\varepsilon_0,\qquad {\rm for\ all}\
t\in\mathbb{R},$$ which combined with Corollary~\ref{co:decay}
implies that for any $\sigma\in\mathbb{R}$ there exists a constant
$C_{H,\sigma}>0$ such that:
 \begin{equation}\label{est:psieh2}
 \max\left\{\|<x>^\sigma\psi_E(t)\|_{H^2},\|<x>^\sigma\frac{\partial\psi_E}{\partial
a_1}(t)\|_{H^2},\|<x>^\sigma\frac{\partial\psi_E}{\partial
a_2}(t)\|_{H^2}\right\}\le C_{H,\sigma}\varepsilon_0,\qquad {\rm
for\ all}\ t\in\mathbb{R}.
 \end{equation}
Consequently, using the continuous imbedding
$H^2(\mathbb{R}^2)\hookrightarrow L^p(\mathbb{R}^2),\ 2\le
p\le\infty$ and $L^2_{\sigma}(\mathbb{R}^2)\hookrightarrow
L^1(\mathbb{R}^2),\ \sigma>1$ we have that for all $1\le p\le\infty$
and all $\sigma\in\mathbb{R},$ there exists the constants
$C_{p,\sigma}$ such that
\begin{equation}\label{est:pPsi}
\sup_{t\in\mathbb{R}}\max\left\{\|\psi_E(t)\|_{L^p_\sigma},\|\frac{\partial\psi_E}{\partial
a_1}(t)\|_{L^p_\sigma},\|\frac{\partial\psi_E}{\partial
a_2}(t)\|_{L^p_\sigma},\|\Psi_1(a(t))\|_{L^p_\sigma},\|\Psi_1(a(t))\|_{L^p_\sigma}\right\}\le
C_{p,\sigma}\varepsilon_0,\end{equation} see Remark \ref{rmk:inv}
for the definitions of $\Psi_j(a),\ j=1,2.$
 In addition,
 since
  $$u\in C(\mathbb{R},H^{1}(\mathbb{R}^2))\cap
  C^1(\mathbb{R},H^{-1}(\mathbb{R}^2)),$$
and $u\mapsto a$ respectively $u\mapsto \eta$ are $C^1,$ see Lemma
\ref{lem:decomp}, we get that $a(t)$ is $C^1$ and $\eta\in
C(\mathbb{R},H^{1})\cap
  C^1(\mathbb{R},H^{-1}).$

The solution is now described by the $C^1$ function
$a:\mathbb{R}\mapsto\mathbb{C}$ and $\eta(t)\in
C(\mathbb{R},H^1)\cap C^1(\mathbb{R},H^{-1}).$ To obtain estimates
for them it is useful to first remove their dominant phase. Consider
the $C^2$ function:
\begin{equation}\label{def:theta}
\theta(t)=\int_0^tE(|a(s)|)ds
\end{equation}
and
 \begin{equation}\label{def:tu}
 \tilde u(t)=e^{i\theta(t)}u(t),
 \end{equation}
then $\tilde u(t)$ satisfies the differential equation:
 \begin{equation}\label{eq:tu}
 i\partial\tilde u(t)=-E(|a(t)|)\tilde u(t)+(-\Delta+V)\tilde u(t)+g(\tilde u(t)),
 \end{equation}
see \eqref{eq:ufull} and \eqref{gsym}. Moreover, like $u(t),$
$\tilde u(t)$ can be decomposed:
 \begin{equation}\label{decomp:tu}
 \tilde u(t)=\underbrace{\tilde a(t)\psi_0+h(\tilde a(t))}_{\tilde\psi_E(t)}+\tilde\eta(t)
 \end{equation}
where
 \begin{equation}\label{def:taeta}
 \tilde a(t)=e^{i\theta(t)}a(t),\qquad
 \tilde\eta(t)=e^{i\theta(t)}\eta(t)\in {\cal H}_{\tilde a(t)}
 \end{equation}
By plugging in \eqref{decomp:tu} into \eqref{eq:tu} we get
 \begin{eqnarray}
 i\frac{\partial\tilde\eta}{\partial
 t}+iD\tilde\psi_E|_{\tilde a} \frac{d\tilde a}{dt}&=&(-\Delta+V-E(|a|)(\tilde\psi_E+\tilde\eta)+g(\tilde\psi_E)+g(\tilde\psi_E+\tilde\eta)-g(\tilde\psi_E)\nonumber\\
 &=&L_{\tilde\psi_E}\tilde\eta+g_2(\tilde\psi_E,\tilde\eta)\nonumber
 \end{eqnarray}
or, equivalently,
 \begin{equation}\label{eq:tudecomp}
 \frac{\partial\tilde\eta}{\partial t}+
 \underbrace{\frac{\partial\tilde\psi_E}{\partial a_1}\frac{d\tilde a_1}{dt}+
 \frac{\partial\tilde\psi_E}{\partial a_2}\frac{d\tilde a_2}{dt}}
 _{\in {\rm span}_\mathbb{R}\{\frac{\partial\tilde\psi_E}{\partial a_1},\frac{\partial\tilde\psi_E}{\partial
 a_2}\}}=\underbrace{-iL_{\tilde\psi_E}\tilde\eta}_{\in {\cal H}_{\tilde a}}-ig_2(\tilde\psi_E,\tilde\eta)
 \end{equation}
where $L_{\tilde\psi_E}$ is defined by \eqref{def:linop}:
 $$L_{\tilde\psi_E}\tilde\eta=
 (-\Delta+V-E(|\tilde a|))\tilde\eta+\frac{d}{d\varepsilon}g(\tilde\psi_E+\varepsilon\tilde\eta)|_{\varepsilon=0}$$
and we used $|a|=|\tilde a|,$ while $g_2$ is defined by:
 \begin{equation}\label{g2}
 g_2(\tilde\psi_E,\tilde\eta)=
 g(\tilde\psi_E+\tilde\eta)-g(\tilde\psi_E)-\frac{d}{d\varepsilon}g(\tilde\psi_E+\varepsilon\tilde\eta)|_{\varepsilon=0}
 \end{equation} and we also used the fact that $\tilde\psi_E$ is a
solution of the eigenvalue problem \eqref{eq:ev}. Note that $g_2$ is
at least quadratic in the second variable, more precisely:

\begin{lemma}\label{le:g2est} There exists a constant $C>0$ such that for all $a,z\in\mathbb{C}$ we have:
 $$|g_2(a,z)|\le
 C(|a|^{\alpha_1}+|a|^{\alpha_2}+|z|^{\alpha_1}+|z|^{\alpha_2})|z|^2$$
\end{lemma}

\smallskip

\noindent{\bf Proof:} From the definition \eqref{g2} of $g_2$ we
have:
 $$
 g_2(a,z)=g(a+z)-g(a)-Dg_{a}[z]=\int_0^1\left(Dg_{a+\tau z}-Dg_{a}\right)[z]d\tau
  =\int_0^1\int_0^1D^2g_{a+s\tau z}[\tau z][z]d\tau ds.$$
Now \eqref{gest} and \eqref{gsym} imply that there exists a constant
$C_1>0$ such that the bilinear form $Dg$ on $\mathbb{C}\times
\mathbb{C}$ satisfies:
 \begin{equation}\label{est:d2g}
 \|D^2g_b\|\le C_1(|b|^{\alpha_1}+|b|^{\alpha_2}),\qquad \forall
 b\in\mathbb{C}.\end{equation} Hence
 $$|g_2(a,z)|\le C_1\left((2\max (|a|,|z|))^{\alpha_1}+(2\max
 (|a|,|z|))^{\alpha_2}\right)\frac{1}{2}|z|^2,$$ which proves the
lemma.\ $\Box$

\smallskip

We now project \eqref{eq:tudecomp} onto the invariant subspaces of
$-iL_{\tilde\psi_E},$ namely ${\rm span}_\mathbb{R}
\{\frac{\partial\tilde\psi_E}{\partial
a_1},\frac{\partial\tilde\psi_E}{\partial
 a_2}\},$ and ${\cal H}_{\tilde a}.$ More
precisely, we evaluate both the left and right hand side of
\eqref{eq:tudecomp} which are functionals in $H^{-1}(\mathbb{R}^2)$
at $\Psi_j=\Psi_j(\tilde a(t)),\ j=1,2,$ see Remark \ref{rmk:inv},
and take the real parts.  We obtain:
 $$\left[\begin{array}{c}\Re\langle\Psi_1,\frac{\partial\tilde\eta}{\partial
 t}\rangle\\ \Re\langle\Psi_2,\frac{\partial\tilde\eta}{\partial
 t}\rangle\end{array}\right]+\frac{d}{dt}\left[\begin{array}{c}\tilde
 a_1\\ \tilde a_2\end{array}\right]=\left[\begin{array}{c}
 g_{21}(\tilde\psi_E,\tilde\eta)\\
 g_{22}(\tilde\psi_E,\tilde\eta)\end{array}\right]$$
where
 \begin{equation}\label{def:g2j}
 g_{2j}(\tilde\psi_E,\tilde\eta)=\Re\langle\Psi_j,-ig_2(\tilde\psi_E,\tilde\eta)\rangle ,\qquad j=1,2.
 \end{equation}
Note that from Lemma~\ref{le:g2est} and H\" older inequality we have
for all $t\in\mathbb{R}:$
 \begin{eqnarray}
 \lefteqn{|g_{2j}(\tilde\psi_E(t),\tilde\eta(t))|\le
 C\int_{\mathbb{R}^2}|\Psi_j(t,x)|\left(|\tilde\psi_E(t,x)|^{\alpha_1}+|\tilde\psi_E(t,x)|^{\alpha_2}
 +|\tilde\eta(t,x)|^{\alpha_1}+|\tilde\eta(t,x)|^{\alpha_2}\right)\
 |\tilde\eta(t,x)|^{2}dx}\label{est:g2j}\\
 &\le & C\left[
 \|\Psi_j(t)\|_{L^{r_0}}\left(\|\tilde\psi_E(t)\|_{L^\infty}^{\alpha_1}+\|\tilde\psi_E(t)\|_{L^\infty}^{\alpha_2}\right)
 \|\tilde\eta(t)\|_{L^{p_2}}^2+\|\Psi_j(t)\|_{L^{r_1}}\|\tilde\eta(t)\|_{L^{p_2}}^{2+\alpha_1}
 +\|\Psi_j(t)\|_{L^{r_2}}\|\tilde\eta(t)\|_{L^{p_2}}^{2+\alpha_2}\right],\nonumber
 \end{eqnarray}
where $r_0^{-1}+(p_2/2)^{-1}=1,\
r_j^{-1}+(p_2/(2+\alpha_j))^{-1}=1,\ j=1,2.$

To calculate $\Re\langle\Psi_j,\frac{\partial\tilde\eta}{\partial
 t}\rangle,\ j=1,2$ we use the fact that $\tilde\eta(t)\in {\cal
 H}_{\tilde a},$ for all $t\in\mathbb{R},$ i.e. $$\Re\langle\Psi_j(\tilde
 a(t)),\tilde\eta(t)\rangle\equiv 0.$$ Differentiating the latter with
respect to $t$ we get:
 $$\Re\langle\Psi_j,\frac{\partial\tilde\eta}{\partial
 t}\rangle=-\Re\langle\frac{\partial\Psi_j}{\partial a_1}\frac{d\tilde a_1}{dt}+
 \frac{\partial\Psi_j}{\partial a_2}\frac{d\tilde a_2}{dt},\tilde\eta\rangle\qquad j=1,2$$
which replaced into above leads to:
 \begin{equation}\label{eq:ta}
 \frac{d}{dt}\left[\begin{array}{c}\tilde
 a_1\\ \tilde
 a_2\end{array}\right]=(\mathbb{I}_{\mathbb{R}^2}-M_{\tilde u})^{-1}\left[\begin{array}{c}
 g_{21}(\tilde\psi_E,\tilde\eta)\\
 g_{22}(\tilde\psi_E,\tilde\eta)\end{array}\right],
 \end{equation}
where the two by two matrix $M_{\tilde u}$ is defined in Remark
\ref{rmk:inv}. In particular
 $$\left[\begin{array}{c}\Re\langle\Psi_1,\frac{\partial\tilde\eta}{\partial
 t}\rangle\\ \Re\langle\Psi_2,\frac{\partial\tilde\eta}{\partial
 t}\rangle\end{array}\right]=-M_{\tilde u}(\mathbb{I}_{\mathbb{R}^2}-M_{\tilde u})^{-1}\left[\begin{array}{c}
 g_{21}(\tilde\psi_E,\tilde\eta)\\
 g_{22}(\tilde\psi_E,\tilde\eta)\end{array}\right],$$
which we use to obtain the component in ${\cal H}_{\tilde
a}=\{\Psi_1(\tilde a),\Psi_2(\tilde a)\}^\perp$ of
\eqref{eq:tudecomp}:
 $$\frac{\partial\tilde\eta}{\partial
 t}=-iL_{\tilde\psi_E}\tilde\eta-ig_2(\tilde\psi_E,\tilde\eta)
 -(\mathbb{I}-M_{\tilde u})^{-1}g_3(\tilde\psi_E,\tilde\eta),$$
where $g_3$ is the projection of $-ig_2$ onto ${\rm
span}_\mathbb{R}\{\frac{\partial\tilde\psi_E}{\partial
 a_1},\frac{\partial\tilde\psi_E}{\partial
 a_2}\}$ relative to the decomposition \eqref{h-1decomp}:
 \begin{equation}\label{def:g3}
 g_3(\tilde\psi_E,\tilde\eta)=g_{21}(\tilde\psi_E,\tilde\eta)\frac{\partial\tilde\psi_E}{\partial
 a_1}+g_{22}(\tilde\psi_E,\tilde\eta)\frac{\partial\tilde\psi_E}{\partial
 a_2},
 \end{equation}
see \eqref{def:g2j} for the definitions of $g_{2j},\ j=1,2,$ and
$\mathbb{I}-M_{\tilde u}$ is the linear operator on the two
dimensional real vector space ${\rm
span}_\mathbb{R}\{\frac{\partial\tilde\psi_E}{\partial
a_1},\frac{\partial\tilde\psi_E}{\partial a_2}\}$ whose matrix
representation relative to the basis
$\{\frac{\partial\tilde\psi_E}{\partial
a_1},\frac{\partial\tilde\psi_E}{\partial a_2}\}$ is
$\mathbb{I}_{\mathbb{R}^2}-M_{\tilde u}.$ It is easier to switch
back to the variable $\eta(t)=e^{-i\theta(t)}\tilde\eta(t)\in {\cal
H}_a:$
 \begin{equation}\label{eq:eta}
 \frac{\partial\eta}{\partial
 t}=-i(-\Delta+V)\eta-iDg_{\psi_E}\eta-ig_2(\psi_E,\eta)
 -(\mathbb{I}-M_u)^{-1}g_3(\psi_E,\eta),
 \end{equation}
where we used the equivariant symmetry \eqref{gsym} and its obvious
consequences for the symmetries of $Dg,\ g_2,\ g_3$ and $M.$ Since
by Lemma \ref{le:pcinv} it is sufficient to get estimates for
$\zeta(t)=P_c\eta (t),$ we now project \eqref{eq:eta} onto the
continuous spectrum of $-\Delta+V:$
 \begin{equation}\label{eq:zeta}
 \frac{\partial\zeta}{\partial
 t}=-i(-\Delta+V)\zeta-iP_cDg_{\psi_E}R_a\zeta-iP_cg_2(\psi_E,R_a\zeta)
 -P_c(\mathbb{I}-M_u)^{-1}g_3(\psi_E,R_a\zeta),
 \end{equation}
where $R_a:{\cal H}_0\mapsto {\cal H}_a$ is the inverse of $P_c$
restricted to ${\cal H}_a,$ see Lemma \ref{le:pcinv}.

Consider the initial value problem for the linear part of
\eqref{eq:zeta}:
 \begin{eqnarray}
 \frac{\partial z}{\partial
 t}&=&-i(-\Delta+V) z-iP_cDg_{\psi_E(t)}R_{a(t)} z(t)\label{eq:zm}\\
 z(s)&=&v\in {\cal H}_0\nonumber
 \end{eqnarray}
and write its solution in terms of a family of operators:
 \begin{equation}\label{def:Omega}
 \Omega(t,s): {\cal
H}_0\mapsto {\cal H}_0,\qquad \Omega(t,s)v=z(t),\ \ t,\
s\in\mathbb{R}.\end{equation} In Section \ref{se:lin} we show that
such a family of operators exists, is uniformly bounded in $t,\ s$
with respect to the $L^2$ norm and it has very similar properties
with the unitary group of operators $e^{-i(-\Delta+V)(t-s)}P_c$
generated by the Schr\" odinger operator $-i(-\Delta+V)P_c.$ In
particular $\Omega(t,s)$ satisfies certain dispersive decay
estimates in weighted $L^2$ spaces and $L^p,\ p>2$ spaces, see
Theorem \ref{th:lin1} and Theorem \ref{th:lin2}. For all these
results to hold we only need to choose $\varepsilon_0$ small enough
such that $\varepsilon_0 C_{H,4\sigma /3}\le \varepsilon_1,$ where
$\sigma >1 $ and $\varepsilon_1>0$ are fixed in Section 4 and the
constant $C_{H,4\sigma /3}$ is the one from \eqref{est:psieh2}.

Using Duhamel formula, the solution $\zeta\in C(\mathbb{R},H^{1}\cap
{\cal H}_0)\cap
  C^1(\mathbb{R},H^{-1}(\mathbb{R}^2)\cap {\cal H}_0)$ of \eqref{eq:zeta}
also satisfies:
\begin{eqnarray}
\zeta(t)&=&\Omega(t,0)\zeta(0)-i\int_0^t\Omega(t,s)P_cg_2(\psi_E(s),R_{a(s)}\zeta(s))ds\nonumber\\
 &&-\int_0^t\Omega(t,s)P_c(\mathbb{I}-M_{u(s)})^{-1}g_3(\psi_E(s),R_{a(s)}\zeta(s))ds.\label{int:zeta}
\end{eqnarray}
Note that the right hand side of \eqref{int:zeta} contains only
terms that are quadratic and higher order in $\zeta,$ see
Lemma~\ref{le:g2est} and \eqref{est:g2j}. As in
\cite{kz:as2d,km:as3d} this is essential in controlling low power
nonlinearities and it is the main difference between our approach
and the existing literature on asymptotic stability of coherent
structures for dispersive nonlinear equations, see \cite[p.
449]{kz:as2d} for a more detailed discussion.

To obtain estimates for $\zeta$ we apply a contraction mapping
argument to the fixed point problem (\ref{int:zeta}) in the
following Banach space. Fix $p_0>2$ such that
 \begin{equation}\label{def:p0}
 p_0>\max\left\{\frac{2}{\alpha_1-1/2},\
 (4+2\alpha_2)\frac{q_0-2}{q_0-(4+2\alpha_2)}\right\},
 \end{equation}
and let
 \begin{equation}\label{def:p2}
 p_2=\frac{p_0q_0}{p_0+q_0-2},
 \end{equation}
and
 \begin{equation}\label{def:p1}
 p_1=\frac{2}{1-\alpha_1+2/p_0},\qquad {\rm if}\ \alpha_1<1,
 \end{equation}
then
\begin{itemize}
 \item[Case I] if $\alpha_1\ge 1,$ or $1/2<\alpha_1<1$ and $p_1\ge
 p_2,$ let:
 $$
 Y=\left\{v\in C(\mathbb{R},L^2\cap L^{p_2})\ :\ \sup_{t\in\mathbb{R}}
 \|v(t)\|_{L^2}<\infty,
 \ \sup_{t\in\mathbb{R}}
 \frac{(1+|t|)^{1-\frac{2}{p_2}}}{[\log(2+|t|)]^{\frac{1-\frac{2}{p_2}}{1-\frac{2}{p_0}}}}\|v(t)\|_{L^{p_2}}<\infty
 \right\};$$
 \item[Case II] if $1/2<\alpha_1<1$ and $p_1<
 p_2,$ let:
 \begin{eqnarray}
 Y&=&\left\{v\in C\left(\mathbb{R},L^2\cap L^{p_1}\cap L^{p_2}\right)\ :\ \sup_{t\in\mathbb{R}}
\|v(t)\|_{L^2}<\infty,\right.\nonumber\\
 &&\sup_{t\in\mathbb{R}}
 \frac{(1+|t|)^{1-\frac{2}{p_1}}}{[\log(2+|t|)]^{\frac{1-\frac{2}{p_1}}{1-\frac{2}{p_0}}}}\|v(t)\|_{L^{p_1}}<\infty,\
 \left.\sup_{t\in\mathbb{R}}
\frac{(1+|t|)^{\alpha_1-\frac{2}{p_0}}}{[\log(2+|t|)]^{\frac{\alpha_1-\frac{2}{p_0}}{1-\frac{2}{p_0}}}}\|v(t)\|_{L^{p_2}}<\infty
 \right\};\nonumber\end{eqnarray}
 \end{itemize}
endowed with the norm $$
 \|v\|_{Y}=\max\left\{\sup_{t\in\mathbb{R}} \|v(t)\|_{L^2},\
\sup_{t\in\mathbb{R}}
\frac{(1+|t|)^{1-\frac{2}{p_2}}}{[\log(2+|t|)]^{\frac{1-\frac{2}{p_2}}{1-\frac{2}{p_0}}}}\|v(t)\|_{L^{p_2}}
 \right\}$$ in Case I, while in Case II
 $$
  \|v\|_{Y}=\max\left\{\sup_{t\in\mathbb{R}}
\|v(t)\|_{L^2},\ \sup_{t\in\mathbb{R}}
\frac{(1+|t|)^{1-\frac{2}{p_1}}}{[\log(2+|t|)]^{\frac{1-\frac{2}{p_1}}{1-\frac{2}{p_0}}}}\|v(t)\|_{L^{p_1}},\
\sup_{t\in\mathbb{R}}
\frac{(1+|t|)^{\alpha_1-\frac{2}{p_0}}}{[\log(2+|t|)]^{\frac{\alpha_1-\frac{2}{p_0}}{1-\frac{2}{p_0}}}}\|v(t)\|_{L^{p_2}}
 \right\}.$$

Consider now the nonlinear operator in (\ref{int:zeta}):
 \begin{equation}\label{def:N}
 N(v)(t)=-i\int_0^t\Omega(t,s)P_cg_2(\psi_E(s),R_{a(s)}v(s))ds
 -\int_0^t\Omega(t,s)P_c(\mathbb{I}-M_{u(s)})^{-1}g_3(\psi_E(s),R_{a(s)}v(s))ds.\end{equation} We
have:

\begin{lemma}\label{lm:lip}
 $N : Y\rightarrow Y$ is well defined, and locally Lipschitz,
 i.e. there exists $\tilde{C}>0$, such that $$\|Nv_1
-Nv_2\|_{Y}\leq\tilde{C}(\|v_1\|_{Y}+\|v_2\|_{Y}+\|v_1\|_{Y}^{1+\alpha_1}+\|v_2\|_{Y}^{1+\alpha_1}+
\|v_1\|_{Y}^{1+\alpha_2}+\|v_2\|_{Y}^{1+\alpha_2})\|v_1 -v_2\|_{Y}.
$$
\end{lemma}

Assuming that the lemma has been proven then we can apply the
contraction principle for \eqref{int:zeta} in a closed ball in the
Banach space $Y$ in the following way. Let
$$v=\Omega(t,0)\zeta(0)$$
then by Theorem \ref{th:lin2}
$$\|v\|_Y\le
\max\{C_2,C_{p_0,p_1},C_{p_0,p_2}\}\|\zeta(0)\|_{L^2\cap L^{q_0'}}$$
where we used the interpolation
$\|\zeta(0)\|_{L^r}\le\|\zeta(0)\|_{L^2\cap L^{q_0'}},\ q_0'\le r\le
2$ with $r=q'$ and $r=p'$ defined in theorem \ref{th:lin2} for
$p=p_j,\ j=1,2.$ Recall that
  $$\zeta(0)=P_c\eta(0)=P_cu_0-h(a(0))=u_0-\langle\psi_0,u_0\rangle\psi_0-h(a(0))$$
where $u_0=u(0)$ is the initial data, see also \eqref{dc}. Hence
 $$\|\zeta(0)\|_{L^2\cap L^{q_0'}}\le\|u_0\|_{L^2\cap
 L^{q_0'}}+\|u_0\|_{L^2}\|\psi_0\|_{L^{q_0'}}+D_1\|u_0\|_{L^2}\le
 D\varepsilon_0$$ where $D_1,\ D>0$ are constants independent on
$u_0$ and the estimate on $h(a(0))$ follows from Proposition
\ref{pr:cm} and $|a(0)|\le 2\|u_0\|_{L^2}$ see Lemma
\ref{lem:decomp}.

Therefore we can choose $\varepsilon_0$ small enough such that
$R=2\|v\|_Y$ satisfies
$$Lip\stackrel{def}{=}2\tilde{C}(R+R^{1+\alpha_1}+
R^{1+\alpha_2})<1.$$ In this case the integral operator given by the
right hand side of \eqref{int:zeta}:
$$K(\zeta)=v+N(\zeta)$$
leaves $B(0,R)={\zeta\in Y:\|\zeta\|_Y\le R}$ invariant and it is a
contraction on it with Lipschitz constant $ Lip$ defined above.
Consequently the equation \eqref{int:zeta} has a unique solution in
$B(0,R)$ and because $\zeta(t)\in C(\mathbb{R},H^1)\hookrightarrow
C(\mathbb{R},L^2,L^{p_1},L^{p_2})$ already verified the equation we
deduce that $\zeta(t)$ is in $B(0,R)$, in particular it satisfies
the estimates \eqref{lpdecay}.

Then $\eta(t)=R_a(t)\zeta(t)$ satisfies \eqref{lpdecay} because of
Lemma \ref{le:pcinv}. Moreover, the system of ODE's \eqref{eq:ta}
has integrable in time right hand side because the matrix has norm
bounded by $2,$ see \eqref{M-bound}, while $g_{2j}$ satisfy
\eqref{est:g2j} where $\tilde\eta(t)$ differs from $\eta(t)$ by only
a phase and the $L^p,\ 1\le p\le\infty$ norms of $\Psi_j(t),\
\psi_E(t)$ are uniformly bounded in time, see \eqref{est:pPsi}.
Consequently $\tilde a_1(t)$ and $\tilde a_2(t)$ converge as
$t\rightarrow\pm\infty,$ and there exist the constants $C, \epsilon
>0$ such that:
 $$\lim_{t\rightarrow\pm\infty}\tilde a(t)=\lim_{t\rightarrow\pm\infty}\tilde
 a_1(t)+i\tilde a_2(t)\stackrel{def}{=} a_{\pm\infty},\qquad
 |\tilde a(\pm t)-a_{\pm\infty}|\le C(1+t)^{-\epsilon},\ {\rm for\
 all}\ t\ge 0.$$

We can now define
 \begin{equation}\label{def:psieinfty}
 \psi_{E_{\pm\infty}}=a_{\pm\infty}\psi_0+h(a_{\pm\infty}),
 \end{equation}
and we have
 \begin{equation}\label{conv:tpsie}\lim_{t\rightarrow\pm\infty}\|\tilde\psi_E(t)-\psi_{E_{\pm\infty}}\|_{H^2\cap
 L^2_\sigma}=0,\ {\rm for}\ \sigma\in\mathbb{R}
 \end{equation}
where we used \eqref{decomp:tu} and the continuity of $h(a),$ see
Proposition \ref{pr:cm}. In addition, since
$E:[-2\delta_1,\delta_1]\mapsto (-\delta,\delta)$ is a $C^1$
function, see Proposition \ref{pr:cm}, the following limits exist
together with the constant $C_1>0$ such that:
 $$\lim_{t\rightarrow\pm\infty}E(|\tilde a(t)|)=E_{\pm\infty},\qquad
 |E(|\tilde a(\pm t)|)-E_{\pm\infty}|\le C_1(1+t)^{-\epsilon} \ {\rm for\
 all}\ t\ge 0.$$
If we define
 \begin{equation}\label{def:ttheta}
 \tilde\theta(t)=\left\{\begin{array}{lr}\frac{1}{t}\int_0^{t}E(|\tilde
 a(s)|)-E_{+\infty}ds & {\rm if}\ t>0\\ 0 & {\rm if}\
 t=0\\ \frac{1}{t}\int_0^{t}E(|\tilde
 a(s)|)-E_{-\infty}ds & {\rm if}\ t<0\end{array}\right.
 \end{equation}
then
 $$\lim_{|t|\rightarrow\infty}\tilde\theta(t)=0$$
and
 \begin{equation}\label{rel:theta1}\theta(t)=\int_0^{t}E(|a(s)|)ds=\left\{\begin{array}{lr}
 t(E_{+\infty}+\tilde\theta(t)) & {\rm if}\ t\ge 0\\ t(E_{+\infty}+\tilde\theta(t)) & {\rm if}\ t< 0
 \end{array}\right.\end{equation}
where we used $|a(t)|=|\tilde a(t)|,$ see \eqref{def:taeta}.

In conclusion, since $\psi_E(t)=e^{i\theta(t)}\tilde\psi_E(t),$ see
\eqref{dc}, \eqref{decomp:tu} and \eqref{def:taeta}, we get from
\eqref{conv:tpsie} and \eqref{rel:theta1} the convergence
\eqref{conv:psie}.

It remains to prove Lemma \ref{lm:lip}:

\smallskip

\noindent{\bf Proof of Lemma \ref{lm:lip}:} It suffices to prove the
estimate:
 \begin{equation}\label{est:N}
 \|Nv_1-Nv_2\|_{Y}\leq\tilde{C}(\|v_1\|_{Y}+\|v_2\|_{Y}+\|v_1\|_{Y}^{1+\alpha_1}+\|v_2\|_{Y}^{1+\alpha_1}+
 \|v_1\|_{Y}^{1+\alpha_2}+\|v_2\|_{Y}^{1+\alpha_2})\|v_1 -v_2\|_{Y},
 \end{equation}
because plugging in $v_2\equiv 0$ and using $N(0)\equiv 0,$ see
\eqref{def:N}, will then imply $N(v_1)\in Y$ whenever $v_1\in Y.$

Note that via interpolation in $L^p$ spaces we have for all $v\in Y$
and any $2\le p\le p_2:$
 \begin{equation}\label{est:lp}
  \|v(t)\|_{L^p}\le \left\{\begin{array}{cl}
     \|v\|_Y\frac{\log^{\frac{1-2/p}{1-2/p_0}}(2+|t|)}{(1+|t|)^{1-2/p}}
     & {\rm if}\ \alpha_1\ge 1\ {\rm or}\ \alpha_1<1\ {\rm and}\ \ p\le \frac{2}{1-\alpha_1+2/p_0}, \\
      & \\
     \|v\|_Y\frac{\log^{\frac{\alpha_1-2/p_0}{1-2/p_0}}(2+|t|)}{(1+|t|)^{\alpha_1-2/p_0}}
     & {\rm if}\ \alpha_1<1\ {\rm and}\ p>
     \frac{2}{1-\alpha_1+2/p_0}.
     \end{array}\right.
 \end{equation}
Now, from \eqref{g2}, we have for any $v_1,\ v_2\in Y:$
 \begin{eqnarray}
 g_2(\psi_E,R_av_1)-g_2(\psi_E,R_av_2)&=&g(\psi_E+R_av_1)-g(\psi_E+R_av_2)-Dg_{\psi_E}[R_a(v_1-v_2)]\nonumber\\
  &=&\int_0^1\left(Dg_{\psi_E+R_a(\tau
  v_1+(1-\tau)v_2)}-Dg_{\psi_E}\right)[R_a(v_1-v_2)]d\tau\nonumber\\
  &=&\int_0^1\int_0^1D^2g_{\psi_E+sR_a(\tau
  v_1+(1-\tau)v_2)}[R_a(\tau v_1+(1-\tau)v_2)][R_a(v_1-v_2)]d\tau
  ds\nonumber\\
  &=&
  A_1(\psi_E,v_1,v_2)+A_2(\psi_E,v_1,v_2)+A_3(\psi_E,v_1,v_2),\label{def:A123}
  \end{eqnarray}
where we consider $\chi_j(t,x),\ j=1,2$ to be the characteristic
function of the set $S_1=\{(t,x)\in\mathbb{R}\times\mathbb{R}^2\ :\
|\psi_E(t,x)|\ge\max (|R_{a(t)}v_1(t,x)|,|R_{a(t)}v_2(t,x)|)\},$
respectively\newline $S_2= \{(t,x)\in\mathbb{R}\times\mathbb{R}^2\
:\ \max (|R_{a(t)}v_1(t,x)|,|R_{a(t)}v_2(t,x)|)\le 1\}$ and
 \begin{eqnarray}
 A_1(\psi_E,v_1,v_2)&=&\int_0^1\int_0^1\chi_1 D^2g_{\psi_E+sR_a(\tau
  v_1+(1-\tau)v_2)}[R_a(\tau v_1+(1-\tau)v_2)][R_a(v_1-v_2)]d\tau
  ds,\nonumber\\
 A_2(\psi_E,v_1,v_2)&=&\int_0^1\int_0^1(1-\chi_1)\chi_2 D^2g_{\psi_E+sR_a(\tau
  v_1+(1-\tau)v_2)}[R_a\tau v_1+(1-\tau)v_2)][R_a(v_1-v_2)]d\tau
  ds,\nonumber\\
 A_3(\psi_E,u_1,u_2)&=&\int_0^1\int_0^1(1-\chi_1)(1-\chi_2) D^2g_{\psi_E+sR_a(\tau
  v_1+(1-\tau)v_2)}[R_a(\tau v_1+(1-\tau)v_2)][R_a(v_1-v_2)]d\tau
  ds.\nonumber
  \end{eqnarray}
Note that there exists a constant $C>0$ such that for any $\psi_E,\
v_1,\ v_2\in Y,$ any $t\in\mathbb{R}$ and almost all
$x\in\mathbb{R}^2$ we have the pointwise estimates:
 \begin{eqnarray}
 |A_1(\psi_E(t,x),v_1(t,x),v_2(t,x))|&\le
 &C\left(2^{\alpha_1}|\psi_E(t,x)|^{\alpha_1}+2^{\alpha_2}|\psi_E(t,x)|^{\alpha_2}\right)(|R_{a(t)}v_1(t,x)|+|R_{a(t)}v_2(t,x)|)\nonumber\\
 &&\times |R_{a(t)}(v_1(t,x)-v_2(t,x))|\nonumber\\
 |A_2(\psi_E(t,x),v_1(t,x),v_2(t,x))|&\le
 &2^{\alpha_1}C\left(|R_{a(t)}v_1(t,x)|^{1+\alpha_1}+|R_{a(t)}v_2(t,x)|^{1+\alpha_1}\right)|R_{a(t)}(v_1(t,x)-v_2(t,x))|\nonumber\\
 |A_3(\psi_E(t,x),v_1(t,x),v_2(t,x))|&\le
 &2^{\alpha_2}C\left(|R_{a(t)}v_1(t,x)|^{1+\alpha_2}+|R_{a(t)}v_2(t,x)|^{1+\alpha_2}\right)|R_{a(t)}(v_1(t,x)-v_2(t,x))|\nonumber
 \end{eqnarray}
where we used \eqref{est:d2g}. Consequently, for any
$\sigma\in\mathbb{R}$  there exists a constant $C_{\sigma}>0$ such
that for any $t\in\mathbb{R}:$
 \begin{eqnarray}
 \|A_1(\psi_E(t),v_1(t),v_2(t))\|_{L^2_\sigma}&\le
 &C\|2^{\alpha_1}|\psi_E(t)|^{\alpha_1}+2^{\alpha_2}|\psi_E(t)|^{\alpha_2}\|_{L^s_\sigma}(\|R_{a(t)}v_1(t)\|_{L^{p_2}}+\|R_{a(t)}v_2(t)\|_{L^{p_2}})\nonumber\\
 &&\times \|R_{a(t)}(v_1(t)-v_2(t))\|_{L^{p_2}}\nonumber\\
 &\le &\frac{C_\sigma
 \log^{a_1}(2+|t|)}{(1+|t|)^{b_1}}
 (\|v_1\|_Y+\|v_2\|_Y)\|v_1-v_2\|_Y\label{est:A1}\\
 \end{eqnarray}
where $\frac{1}{s}+\frac{2}{p_2}=\frac{1}{2},$ and, for $\Psi_j,\
j=1,2$ defined in Remark~\ref{rmk:inv}:
 \begin{eqnarray}
 |\Re\langle\Psi_j(a(t)),-iA_1(\psi_E(t),v_1(t),v_2(t))\rangle|&\le &\|\Psi_j(a(t))\|_{L^2_{-\sigma}}\|A_1(\psi_E(t),v_1(t),v_2(t))\|_{L^2_\sigma}\nonumber\\
 &\le & C_{2,-\sigma}\frac{C_\sigma
 \log^{a_1}(2+|t|)}{(1+|t|)^{b_1}}
 (\|v_1\|_Y+\|v_2\|_Y)\|v_1-v_2\|_Y\label{est:g2ja1}
 \end{eqnarray} where
 \begin{equation}\label{def:ab1}
 b_1=\left\{\begin{array}{lr} 2-\frac{4}{p_2} & {\rm in\ Case\ I},\\
                  2\alpha_1-\frac{4}{p_0} & {\rm in\ Case\
                  II},\end{array}\right.\qquad a_1=\left\{\begin{array}{lr} 2\frac{1-2/p_2}{1-2/p_0} & {\rm in\ Case\ I},\\
                  2\frac{\alpha_1-2/p_0}{1-2/p_0} & {\rm in\ Case\
                  II},\end{array}\right.
 \end{equation}
see the definition of the Banach space $Y,$ and we used H\" older
inequality together with \eqref{est:pPsi} and Lemma \ref{le:pcinv}.

Similarly, for any $1\le r'\le 2$ we have $(2+\alpha_1)r'\le
(2+\alpha_2)r'\le p_2,$ hence the above pointwise estimates and
\eqref{est:lp} imply that there exists a constant $C_{r'}>0$ such
that for any $t\in\mathbb{R}:$
 \begin{eqnarray}
 \|A_2(\psi_E(t),v_1(t),v_2(t))\|_{L^{r'}}&\le &
 2^{\alpha_1}C\|\ |R_{a(t)}v_1(t)|^{1+\alpha_1}+|R_{a(t)}v_2(t)|^{1+\alpha_1}\|_{L^{\frac{(2+\alpha_1)r'}{1+\alpha_1}}}\|R_{a(t)}(v_1(t)-v_2(t))\|_{L^{(2+\alpha_1)r'}}\nonumber\\
 &\le & \frac{C_{r'}
 \log^{a_2(r')}(2+|t|)}{(1+|t|)^{b_2(r')}}
 (\|v_1\|_Y^{1+\alpha_1}+\|v_2\|_Y^{1+\alpha_1})\|v_1-v_2\|_Y,\label{est:A2}
 \end{eqnarray}
where
 \begin{equation}\label{def:ab2}
 \begin{array}{lll} b_2(r')=\alpha_1+\frac{2}{r}, & a_2(r')=\frac{\alpha_1+2/r}{1-2/p_0}, & {\rm if}\ \alpha_1\ge 1\ {\rm or}\ \alpha_1<1\ {\rm and}\ \ (2+\alpha_1)r'\le p_1,\\
                  b_2(r')=(2+\alpha_1)(\alpha_1-\frac{2}{p_0}), & a_2(r')=(2+\alpha_1)\frac{\alpha_1-2/p_0}{1-2/p_0}, & {\rm if}\ \alpha_1<1\ {\rm and}\ \ (2+\alpha_1)r'>
                  p_1,
                  \end{array}
 \end{equation} with $1/r+1/r'=1,$ and
 \begin{eqnarray}
 \|A_3(\psi_E(t),v_1(t),v_2(t))\|_{L^{r'}}&\le
 &2^{\alpha_2}C\|\ |R_{a(t)}v_1(t)|^{1+\alpha_2}+|R_{a(t)}v_2(t)|^{1+\alpha_2}\|_{L^{\frac{(2+\alpha_2)r'}{1+\alpha_2}}}\|R_{a(t)}(v_1(t)-v_2(t))\|_{L^{(2+\alpha_2)r'}}\nonumber\\
 &\le &\frac{C_{r'}
 \log^{a_3(r')}(2+|t|)}{(1+|t|)^{b_3(r')}}
 (\|v_1\|_Y^{1+\alpha_2}+\|v_2\|_Y^{1+\alpha_2})\|v_1-v_2\|_Y,\label{est:A3}
 \end{eqnarray}
where
 \begin{equation}\label{def:ab3}
 \begin{array}{lll} b_3(r')=\alpha_2+\frac{2}{r}, & a_3(r')=\frac{\alpha_2+2/r}{1-2/p_0}, & {\rm if}\ \alpha_1\ge 1\ {\rm or}\ \alpha_1<1\ {\rm and}\ \ (2+\alpha_2)r'\le p_1,\\
                  b_3(r')=(2+\alpha_2)(\alpha_1-\frac{2}{p_0}), & a_3(r')=(2+\alpha_2)\frac{\alpha_1-2/p_0}{1-2/p_0}, & {\rm if}\ \alpha_1<1\ {\rm and}\ \ (2+\alpha_2)r'>
                  p_1.
                  \end{array}
 \end{equation} Moreover, using Cauchy-Schwartz inequality and \eqref{est:pPsi} we have:
 \begin{eqnarray}
 |\Re\langle\Psi_j(a(t)),-iA_2(\psi_E(t),v_1(t),v_2(t))\rangle|&\le &\|\Psi_j(a(t))\|_{L^2}\|A_2(\psi_E(t),v_1(t),v_2(t))\|_{L^{2}}\nonumber\\
 &\le & C_{2,0}\frac{C_{2}
 \log^{a_2(2)}(2+|t|)}{(1+|t|)^{b_2(2)}}
 (\|v_1\|_Y^{1+\alpha_1}+\|v_2\|_Y^{1+\alpha_1})\|v_1-v_2\|_Y,\label{est:g2ja2}
 \end{eqnarray}
and
 \begin{equation}\label{est:g2ja3}
 |\Re\langle\Psi_j(a(t)),-iA_3(\psi_E(t),v_1(t),v_2(t))\rangle|\le  C_{2,0}\frac{C_{2}
 \log^{a_3(2)}(2+|t|)}{(1+|t|)^{b_3(2)}}
 (\|v_1\|_Y^{1+\alpha_2}+\|v_2\|_Y^{1+\alpha_2})\|v_1-v_2\|_Y.
 \end{equation}

 Now, from \eqref{def:g3} and \eqref{def:g2j} we have
 \begin{eqnarray}
 \lefteqn{g_3(\psi_E,R_av_1)-g_3(\psi_E,R_av_2)}\nonumber\\
 &=&\Re\langle\Psi_1(a),-i(g_2(\psi_E,R_av_1)-g_2(\psi_E,R_av_2))\rangle
 \frac{\partial\psi_E}{\partial a_1}+\Re\langle\Psi_2(a),-i(g_2(\psi_E,R_av_1)-g_2(\psi_E,R_av_2))\rangle
 \frac{\partial\psi_E}{\partial a_2}\nonumber\\
 &=&\Re\langle\Psi_1(a),-i(A_1+A_2+A_3)(\psi_E,v_1,v_2)\rangle
 \frac{\partial\psi_E}{\partial a_1}+\Re\langle\Psi_2(a),-i(A_1+A_2+A_3)(\psi_E,v_1,v_2))\rangle
 \frac{\partial\psi_E}{\partial a_2}.\nonumber
 \end{eqnarray}
Consequently, for
 \begin{equation}\label{def:A4}
 A_4(\psi_E,v_1,v_2)\stackrel{def}{=}(\mathbb{I}-M_u)^{-1}(g_3(\psi_E,R_av_1)-g_3(\psi_E,R_av_2))
 \end{equation}
we have that for any $\sigma\in\mathbb{R}$ there exists a constant
$C_\sigma>0$ such that:
 \begin{eqnarray}
 \lefteqn{\|A_4(\psi_E(t),v_1(t),v_2(t))\|_{L^2_\sigma}\le \max\left\{\|\frac{\partial\psi_E}{\partial
 a_1}(t)\|_{L^2_\sigma},\|\frac{\partial\psi_E}{\partial
 a_2}(t)\|_{L^2_\sigma}\right\}\sqrt{2}\|(\mathbb{I}-M_{u(t)})^{-1}\|_{\mathbb{R}^2\mapsto\mathbb{R}^2}}\nonumber\\
 &&\times\sqrt{|\Re\langle\Psi_1(a(t)),-i(A_1+A_2+A_3)(t)\rangle|^2+|\Re\langle\Psi_2(a(t)),-i(A_1+A_2+A_3)(t)\rangle|^2}\nonumber\\
 &\le &\frac{C_\sigma
 \log^{a_4}(2+|t|)}{(1+|t|)^{b_4}}
 (\|v_1\|_{Y}+\|v_2\|_{Y}+\|v_1\|_{Y}^{1+\alpha_1}+\|v_2\|_{Y}^{1+\alpha_1}+
 \|v_1\|_{Y}^{1+\alpha_2}+\|v_2\|_{Y}^{1+\alpha_2})\|v_1-v_2\|_Y\label{est:A4}
 \end{eqnarray}
where
 \begin{equation}\label{def:ab4}
 b_4=\min\{b_1,b_2(2),b_3(2)\},\qquad a_4=\max\{a_1,a_2(2),a_3(3)\},
 \end{equation}
and we used \eqref{est:pPsi}, \eqref{M-bound}, \eqref{est:g2ja1},
\eqref{est:g2ja2}, and \eqref{est:g2ja3}.

We are now ready to prove the Lipschitz estimate for the nonlinear
operator $N,$ \eqref{est:N}. From its definition \eqref{def:N} and
\eqref{def:A123}, \eqref{def:A4} we have for any $v_1,v_2\in Y,$ any
$2\le p\le p_2,$ and a fixed $\sigma>1:$
 \begin{eqnarray}
 \lefteqn{\|N(v_1)(t)-N(v_2)(t)\|_{L^p}=\left\|\int_0^t\Omega(t,s)P_c(-iA_1-iA_2-iA_3-A_4)(\psi_E(s),v_1(s),v_2(s))ds\right\|_{L^p}}\nonumber\\
 &\le &\int_0^t\|\Omega(t,s)\|_{L^2_\sigma\mapsto
 L^p}\left(\|A_1(\psi_E(s),v_1(s),v_2(s))\|_{L^2_\sigma}+\|A_4(\psi_E(s),v_1(s),v_2(s))\|_{L^2_\sigma}\right)ds\nonumber\\
 &+&\int_0^{|t|}\|\Omega(t,s)\|_{L^{q'}\cap L^{p'}\mapsto
 L^p}\left(\|A_2(\psi_E(s),v_1(s),v_2(s))\|_{L^{q'}\cap L^{p'}}+\|A_3(\psi_E(s),v_1(s),v_2(s))\|_{L^{q'}\cap L^{p'}}\right)ds.\nonumber
 \end{eqnarray}
where
 \begin{equation}\label{def:p'q'}
 1/p'+1/p=1,\ q'=p'(p_0-2)/(p_0-p'),\ 1/q+1/q'=1.\end{equation} From
Theorem~\ref{th:lin1} and estimates \eqref{est:A1}, \eqref{est:A4}
we get:
 \begin{eqnarray}
 \lefteqn{\int_0^{|t|}\|\Omega(t,s)\|_{L^2_\sigma\mapsto
 L^p}\left(\|A_1(\psi_E(s),v_1(s),v_2(s))\|_{L^2_\sigma}+\|A_4(\psi_E(s),v_1(s),v_2(s))\|_{L^2_\sigma}\right)ds}\nonumber\\
 &\le &(\|v_1\|_{Y}+\|v_2\|_{Y}+\|v_1\|_{Y}^{1+\alpha_1}+\|v_2\|_{Y}^{1+\alpha_1}+
 \|v_1\|_{Y}^{1+\alpha_2}+\|v_2\|_{Y}^{1+\alpha_2})\|v_1-v_2\|_Y\nonumber\\
 &&\times\int_0^t\frac{C_p}{|t-s|^{1-2/p}}\left[\frac{C_\sigma \log^{a_1}(2+|s|)}{(1+|s|)^{b_1}}
 +\frac{C_\sigma
 \log^{a_4}(2+|s|)}{(1+|s|)^{b_4}}\right]ds\nonumber
 \end{eqnarray}
while from Theorem~\ref{th:lin2} and estimates \eqref{est:A2},
\eqref{est:A3} we get:
 \begin{eqnarray}
 \lefteqn{\int_0^{|t|}\|\Omega(t,s)\|_{L^{q'}\cap L^{p'}\mapsto
 L^p}\|A_2(\psi_E(s),v_1(s),v_2(s))\|_{L^{q'}\cap L^{p'}}ds\le
 (\|v_1\|_{Y}^{1+\alpha_1}+\|v_2\|_{Y}^{1+\alpha_1})\|v_1-v_2\|_Y}\nonumber\\
 &&\times\int_0^t\frac{C_{p_0,p}\log^{\frac{1-2/p}{1-2/p_0}}(2+|t-s|)}{|t-s|^{1-2/p}}
 \max\left\{\frac{C_{q'} \log^{a_2(q')}(2+|s|)}{(1+|s|)^{b_2(q')}},\frac{C_{p'}
 \log^{a_2(p')}(2+|s|)}{(1+|s|)^{b_2(p')}}\right\}ds\nonumber
 \end{eqnarray}
and
 \begin{eqnarray}
 \lefteqn{\int_0^{|t|}\|\Omega(t,s)\|_{L^{q'}\cap L^{p'}\mapsto
 L^p}\|A_3(\psi_E(s),v_1(s),v_2(s))\|_{L^{q'}\cap L^{p'}}ds \le
 (\|v_1\|_{Y}^{1+\alpha_1}+\|v_2\|_{Y}^{1+\alpha_1})\|v_1-v_2\|_Y}\nonumber\\
 &&\times\int_0^t\frac{C_{p_0,p}\log^{\frac{1-2/p}{1-2/p_0}}(2+|t-s|)}{|t-s|^{1-2/p}}
 \max\left\{\frac{C_{q'} \log^{a_3(q')}(2+|s|)}{(1+|s|)^{b_3(q')}},\frac{C_{p'}
 \log^{a_3(p')}(2+|s|)}{(1+|s|)^{b_3(p')}}\right\}ds.\nonumber
 \end{eqnarray}

In Case I, i.e. $\alpha_1\ge 1,$ or $1/2<\alpha_1<1$ and $p_1\ge
 p_2,$  since $\alpha_2\ge\alpha_1$
and $p_2\ge 4+2\alpha_2>4,$ we have from \eqref{def:ab1},
\eqref{def:ab2}, \eqref{def:ab3} and \eqref{def:ab4} for
$r'\in\{q',p',2\}$ and $1/r+1/r'=1:$
$$b_1=2-\frac{4}{p_2}>1,\ b_2(r')=\alpha_1+\frac{2}{r}>1,\
b_3(r')=\alpha_2+\frac{2}{q}>1,\ b_4=\min\{b_1,b_2(2),b_3(2)\}>1.$$
We now use the following known convolution estimate:
 \begin{equation}\label{est:conv}
 \int_0^{|t|}\frac{\log^a(2+|t-s|)}{|t-s|^b}\frac{\log^c(2+|s|)}{(1+|s|)^d}ds\le
 C(a,b,c,d)\frac{\log^a(2+|t|)}{(1+|t|)^b},\qquad {\rm for}\ d>1,\
 b<1,\end{equation}
to bound the integral terms above and obtain for all $2\le p\le
p_2:$
 \begin{eqnarray}
 \lefteqn{\|N(v_1)(t)-N(v_2)(t)\|_{L^p}\le
 C_p\frac{\log^{\frac{1-2/p}{1-2/p_0}}(2+|t|)}{(1+|t|)^{1-2/p}}}\nonumber\\
 &&\times (\|v_1\|_{Y}+\|v_2\|_{Y}+\|v_1\|_{Y}^{1+\alpha_1}+\|v_2\|_{Y}^{1+\alpha_1}+
 \|v_1\|_{Y}^{1+\alpha_2}+\|v_2\|_{Y}^{1+\alpha_2})\|v_1-v_2\|_Y\label{est:Np1}\\
 \end{eqnarray}
which, upon moving the time dependent terms to the left hand side
and taking supremum over $t\in\mathbb{R}$ when $p\in\{2,p_2\},$
leads to \eqref{est:N} for $\tilde C=\max\{C_2,C_{p_2}\}.$

In Case II, i.e. $1/2<\alpha_1<1$ and $p_1<
 p_2,$  we have from \eqref{def:ab1}
$b_1=2(\alpha_1-\frac{2}{p_0})>1$ because $p_0>2/(\alpha_1-1/2),$
see \eqref{def:p0}. From \eqref{def:ab2}, under the restriction
$2\le p\le p_1,$ with $p',\ q',\ q$ defined by \eqref{def:p'q'}, we
have either:
 $$b_2(p')>b_2(q')=\alpha_1+2/q>1,$$ or
 $$b_2(p')=b_2(q')=(2+\alpha_1)(\alpha_1-2/p_0)>(2+\alpha_1)/2>1.$$
Since $\alpha_2\ge \alpha_1$ implies $b_3(\cdot)\ge b_2(\cdot)$ we
deduce that, under the restriction $2\le p\le p_1,$ we also have
 $$b_3(p')\ge b_3(q')\ge b_2(q')>1,$$
and
 $$b_4=\min\{b_1,b_2(2),b_3(2)\}>1.$$ We can again apply
\eqref{est:conv} to the above integral terms and get for $2\le p\le
p_1$ the estimate \eqref{est:Np1}. For $p>p_1$ one can show that
$(2+\alpha_1)q'<p_1$ hence $b_2(q')=\alpha_1+2/q,$ and, in the
particular case of $p=p_2,$ we get
 $$b_2(q'_2)=\alpha_1+2/q_2<1,$$
where $q'_2,\ q_2$ are given by \eqref{def:p'q'}. We now have from
convolution estimates:
 $$\int_0^{|t|}\frac{\log^{\frac{1-2/p_2}{1-2/p_0}}(2+|t-s|)}{|t-s|^{1-2/p_2}}\frac{\log^{a_2(q'_2)}(2+|s|)}{(1+|s|)^{b_2(q'_2)}}ds
 \le
 C(p_2)\frac{\log^{\frac{1-2/p_2}{1-2/p_0}+a_2(q'_2)}(2+|t|)}{(1+|t|)^{\alpha_1+2/q_2-2/p_2}} \le
 \tilde C(p_2)\frac{\log^{\frac{\alpha_1-2/p_0}{1-2/p_0}}(2+|t|)}{(1+|t|)^{\alpha_1-2/p_0}},$$
where we used \eqref{def:p'q'} and $p_2\le p_0$ to obtain:
 $$\frac{2}{p_2}-\frac{2}{q_2}=\frac{2}{p_0}\left(\frac{1-2/p_2}{1-2/p_0}\right)\le
 \frac{2}{p_0}.$$ Since $b_2(p'_2)>b_2(q'_2)$ and $b_3(p'_2)\ge b_3(q'_2)\ge
 b_2(q'_2)$ we deduce
 \begin{eqnarray}
 \lefteqn{\|N(v_1)(t)-N(v_2)(t)\|_{L^{p_2}}\le
 \tilde C_{p_2}\frac{\log^{\frac{\alpha_1-2/p_0}{1-2/p_0}}(2+|t|)}{(1+|t|)^{\alpha_1-2/p_0}}}\nonumber\\
 &&\times (\|v_1\|_{Y}+\|v_2\|_{Y}+\|v_1\|_{Y}^{1+\alpha_1}+\|v_2\|_{Y}^{1+\alpha_1}+
 \|v_1\|_{Y}^{1+\alpha_2}+\|v_2\|_{Y}^{1+\alpha_2})\|v_1-v_2\|_Y\nonumber\\
 \end{eqnarray}
which, combined with \eqref{est:Np1} for $p\in\{2,p_1\},$ after
moving the time dependent terms on the left hand side and taking
supremum over $t\in\mathbb{R},$ gives \eqref{est:N} in the Case II
with $\tilde C=\max\{C_2,C_{p_1},\tilde C_{p_2}\}.$

This finishes the proof of Lemma~\ref{lm:lip} and of
Theorem~\ref{th:main}. $\Box$

\section{Linear Estimates}\label{se:lin}

\par Consider the linear Schr\" odinger equation with a potential in
two space dimensions: \[
  \begin{cases}
  i\frac{\partial u}{\partial t}=(-\Delta+V(x))u\\
  u(0)=u_0.
  \end{cases}
 \]
It is known that if $V$ satisfies hypothesis (H1)(i) and (ii) then
the radiative part of the solution, i.e. its projection onto the
continuous spectrum of $H=-\Delta +V,$ satisfies the estimates:
 \begin{equation}\label{Murata}
 \|e^{-iHt}P_c u_0\|_{L^2_{-\sigma}}\le C_M
 \frac{1}{(1+|t|)\log^2(2+|t|)}\|u_0\|_{L^2_\sigma},\qquad t\in\mathbb{R},
 \end{equation}
for any $\sigma >1$ and some constant $C_M>0$ depending only on
$\sigma$ see \cite[Theorem 7.6 and Example 7.8]{mm:ae}, and
 \begin{equation}\label{est:Lp}
 \|e^{-iHt}P_c u_0\|_{L^p}\le  \frac{C_p}{|t|^{1-2/p}}\|u_0\|_{L^{p'}}
 \end{equation}
for some constant $C_p>0$ depending only on $p\ge 2$ and $p'$ given
by $p'^{-1}+p^{-1}=1.$ The case $p=\infty$ in (\ref{est:Lp}) is
proven in \cite{ws:de2}. The conservation of the $L^2$ norm, see
\cite[Corollary 4.3.3]{caz:bk}, gives the $p=2$ case:
$$\|e^{-iHt}P_c u_0\|_{L^2}=
\|u_0\|_{L^{2}}.$$ The general result (\ref{est:Lp}) follows from
Riesz-Thorin interpolation.

\par We would like to extend these estimates to the linearized dynamics
around the center manifold. In other words we consider the linear
equation \eqref{eq:zm}, with initial data at time
  $s:$\begin{eqnarray}
 \frac{\partial z}{\partial
 t}&=&-i(-\Delta+V) z-iP_cDg_{\psi_E(t)}R_{a(t)} z(t)\nonumber\\
 z(s)&=&v\in {\cal H}_0\nonumber
 \end{eqnarray}
Note that this is a nonautonomous problem as the bound state
$\psi_E$ around which we linearize may change with time.

  \par By Duhamel's principle we have:
 \begin{eqnarray}
 z(t)=e^{-iH(t-s)}P_c v-i\int_s^t e^{-iH(t-\tau)}
P_cDg_{\psi_E(\tau)}R_{a(\tau)} z(\tau)d\tau
 \label{rel:Duhamellin}
 \end{eqnarray}

 \par As in (\ref{def:Omega}) we denote
 \begin{equation}\label{def:Omega1}
 \Omega(t,s)v\stackrel{def}{=}z(t).
 \end{equation}
Relying on the fact that $\psi_E(t)$ is small and localized
uniformly in $t\in\mathbb{R},$ we have shown in \cite[Section
4]{kz:as2d} for the particular case of cubic nonlinearity,
$g(s)=s^3,\ s\in\mathbb{R},$ that estimates of type
(\ref{Murata})-(\ref{est:Lp}) can be extended to the operator
$\Omega(t,s).$ Due to \eqref{est:dg} which implies for $\sigma\ge 0$
and $1\le p'\le 2:$
 \begin{eqnarray}
 \|Dg_{\psi_E}R_{a} z\|_{L^2_\sigma}&\le &
 C\left(\|\psi_E\|_{L^\infty_{2\sigma/(1+\alpha_1)}}^{1+\alpha_1}+\|\psi_E\|_{L^\infty_{2\sigma/(1+\alpha_2)}}^{1+\alpha_2}\right)\
 C_{-\sigma}\|z\|_{L^2_{-\sigma}}\label{est:dgls}\\
 \|Dg_{\psi_E}R_{a} z\|_{L^{p'}}&\le &
 C\left(\|\psi_E\|_{L^{(1+\alpha_1)q}_{\sigma/(1+\alpha_1)}}^{1+\alpha_1}+\|\psi_E\|_{L^{(1+\alpha_2)q}_{\sigma/(1+\alpha_2)}}^{1+\alpha_2}\right)\
 C_{-\sigma}\|z\|_{L^2_{-\sigma}},\quad
 \frac{1}{p'}=\frac{1}{q}+\frac{1}{2}\label{est:dglps}\\
 \|Dg_{\psi_E}R_{a} z\|_{L^{p'}}&\le &
 C\left(\|\psi_E\|_{L^{(1+\alpha_1)q}}^{1+\alpha_1}+\|\psi_E\|_{L^{(1+\alpha_2)q}}^{1+\alpha_2}\right)\
 C_{r}\|z\|_{L^r},\quad
 \frac{1}{p'}=\frac{1}{q}+\frac{1}{r}\label{est:dglp}
 \end{eqnarray}
see also Lemma \ref{le:pcinv}, we can use, with obvious
modifications, the arguments in \cite[Section 4]{kz:as2d} to show
that:

\begin{theorem} \label{th:lin1} Fix $\sigma >1.$ There exists $\varepsilon_1>0$ such that
if $\|<x>^{4\sigma /3}\psi_E(t)\|_{H^2}<\varepsilon_1$ for all
$t\in\mathbb{R},$ then there exist constants $C,\ C_p>0$ with the
property that for any $t,\ s\in\mathbb{R}$ the following hold:
 \begin{eqnarray}
 \|\Omega(t,s)\|_{L^2_\sigma\mapsto L^2_{-\sigma}}&\le &
\frac{C}{(1+|t-s|)\log^2(2+|t-s|)},\nonumber\\
 \|\Omega(t,s)\|_{L^{p'}\mapsto L^2_{-\sigma}}&\le
&\frac{C_p}{|t-s|^{1-\frac{2}{p}}},\ {\rm for\ any}\ 2\le p<\infty\
{\rm where}\ p'^{-1}+p^{-1}=1,\nonumber\\
\|\Omega(t,s)\|_{L^2_\sigma\mapsto L^p}&\le &
\frac{C_p}{|t-s|^{1-\frac{2}{p}}},\ {\rm for\ any}\ p\ge 2
 \end{eqnarray}
\end{theorem}
and, for:
\begin{equation}\label{def:T}
T(t,s)=\Omega(t,s)-e^{-iH(t-s)}P_c,
\end{equation}
\begin{lemma}\label{le:T} Assume that
$\|<x>^{4\sigma /3}\psi_E(t)\|_{H^2}<\varepsilon_1,\
t\in\mathbb{R},$ where $\varepsilon_1$ is the one used in
Theorem~\ref{th:lin1}. Then for each $1< q'\le 2$ and $2<p<\infty$
there exist the constants $C_{q'},\ C_{p,q'}>0$ such that for all
$t,\ s\in\mathbb{R}$ we have:
\begin{eqnarray}
 \|T(t,s)\|_{L^1\cap L^{q'}\mapsto L^2_{-\sigma}}&\le &
 \frac{C_{q'}}{1+|t-s|},\nonumber\\
 \|T(t,s)\|_{L^1\cap L^{q'}\mapsto
L^p}&\le &\frac{C_{p,q'}\log(2+|t-s|)}
{(1+|t-s|)^{1-\frac{2}{p}}}.\nonumber
 \end{eqnarray}
\end{lemma}
Note that according to the proofs in \cite[Section 4]{kz:as2d}
$C_{q'}\rightarrow\infty$ as $q'\rightarrow 1$ and
$C_{p,q'}\rightarrow\infty$ as $q'\rightarrow 1$ or
$p\rightarrow\infty.$ These could be prevented and an estimate of
the type
 \begin{equation}\label{est:Timp}
 \|T(t,s)\|_{L^1\mapsto
L^\infty}\le \frac{C\log(2+|t-s|)} {1+|t-s|}
\end{equation}
can be obtained by avoiding the singularity of
$\|e^{-iHt}\|_{L^1\mapsto L^\infty}\sim t^{-1}$ at $t=0$ via a
generalized Fourier multiplier technique developed in \cite[Appendix
and Section 4]{km:as3d}. We choose not to use it here because it
requires stronger restrictions on the potential $V(x)$ like its
Fourier transform should be in $L^1$ while its gradient should be in
$L^p,$ for some $p\ge 2,$ and convergent to zero as
$|x|\rightarrow\infty.$

We now present an improved $L^2$ estimate for the family of
operators $T(t,s):$
\begin{lemma}\label{le:T2} Assume that
$\|<x>^{4\sigma /3}\psi_E(t)\|_{H^2}<\varepsilon_1,\
t\in\mathbb{R},$ where $\varepsilon_1$ is the one used in
Theorem~\ref{th:lin1}. Then there exists the constants $C_{2}>0$
such that for all $t,\ s\in\mathbb{R}$ we have: $$
 \|T(t,s)\|_{L^2\mapsto L^2}\le C_2$$
\end{lemma}

\smallskip\par {\bf Proof:} We are going to use a Kato type
smoothing estimate:
\begin{equation}\label{est:kato}
\|<x>^{-\sigma}e^{-iH t}P_c f(x)\|_{L^2_t(\mathbb{R},L^2_x)}\le
C_{K}\|f\|_{L^2},
\end{equation}
see for example \cite{kn:RS4}. We claim that the previous estimate
still holds if we replace $e^{-iH(t-s)}P_c$ by $\Omega(t,s)$,
namely, there exists a constant $\tilde C_K>0$ such that for any
$s\in\mathbb{R}:$
\begin{equation}\label{est:katoo}
\|<x>^{-\sigma}\Omega(\cdot,s)f\|_{L^2_t(\mathbb{R},L^2_x)}\le
\tilde C_{K} \|f\|_{L^2}.
\end{equation}
Indeed, from \eqref{def:Omega1} and \eqref{rel:Duhamellin}, we have
$$<x>^{-\sigma}\Omega(t,s)v=<x>^{-\sigma} e^{-H(t-s)}P_c v+\int_s^t
<x>^{-\sigma}
e^{-iH(t-\tau)}P_cDg_{\psi_E(\tau)}[R_a(\tau)\Omega(\tau,s)v]d\tau$$
and using \eqref{est:dgls}:
 \begin{eqnarray}
\|\Omega(t,s)v\|_{L^2_{-\sigma}}&\le &\|e^{-H(t-s)}P_c
v\|_{L^2_{-\sigma}} +\int_s^t
\|e^{-iH(t-\tau)}P_c\|_{L^2_\sigma\mapsto
L^2_{-\sigma}}\|Dg_{\psi_E(\tau)}\Omega(\tau,s)v(s)\|_{L^2_\sigma}d\tau \nonumber\\
&\le
&\|e^{-iH(t-s)}v\|_{L^2_{-\sigma}}
+C\sup_{\tau\in\mathbb{R}}\left(\|\psi_E(\tau)\|_{L^\infty_{2\sigma/(1+\alpha_1)}}^{1+\alpha_1}+\|\psi_E(\tau)\|_{L^\infty_{2\sigma/(1+\alpha_2)}}^{1+\alpha_2}\right)\nonumber\\
&\times &\int_\mathbb{R}
\frac{\|\Omega(\tau,s)v(s)\|_{L^2_{-\sigma}}}{
(1+|t-\tau|)\log^2(2+|t-\tau|)}d\tau .\nonumber
 \end{eqnarray}
By Young inequality:
$\|f*g\|_{L^2(\mathbb{R})}\le\|f\|_{L^1(\mathbb{R})}\|g\|_{L^2(\mathbb{R})}$
and \eqref{est:kato} we get
$$\|\Omega(\cdot,s)v\|_{L^2(\mathbb{R},L^2_{-\sigma})}\le
C_K\|v\|_{L^2_x}+C\varepsilon_1
\|\Omega(\cdot,s)v\|_{L^2(\mathbb{R},L^2_{-\sigma})}$$ which implies
\eqref{est:katoo}.

Finally we turn to the estimate in $L^2_x$ for $T(t,s):$
\begin{eqnarray}
\lefteqn{\|T(t,s)v\|_{L^2_x}^2=}\nonumber\\&=& \langle \int_s^t
e^{-iH(t-\tau)}P_cDg_{\psi_E}[R_a\Omega(\tau,s)v]d\tau,
\int_s^t e^{-iH(t-\tau ')} P_cDg_{\psi_E}[R_a\Omega(\tau ',s)v]d\tau '\rangle\nonumber\\
&=&\int_s^t\int_s^t d\tau d\tau' \langle Dg_{\psi_E}[R_a\Omega(\tau,s)v],
e^{-iH(\tau-\tau')}P_cDg_{\psi_E}[R_a\Omega(\tau ',s)v]\rangle\nonumber\\
&\le &C\sup_{\tau\in\mathbb{R}}
\left(\|\psi_E(\tau)\|_{L^\infty_{2\sigma/(1+\alpha_1)}}^{1+\alpha_1}+\|\psi_E(\tau)\|_{L^\infty_{2\sigma/(1+\alpha_2)}}^{1+\alpha_2}\right)^2\nonumber\\
&\times &\int_s^t\int_s^t d\tau d\tau'
\underbrace{\|\Omega(\tau,s)v\|_{L^2_{-\sigma}}}_{\in
L^2(\mathbb{R})}
\underbrace{\|e^{-iH(\tau-\tau')}P_c\|_{L^2_{-\sigma}\mapsto
L^2_{-\sigma}}}_{\in L^1(\mathbb{R})} \underbrace{\|\Omega(\tau
',s)v\|_{L^2_{-\sigma}}}_{\in L^2(\mathbb{R})}. \nonumber
\end{eqnarray}
Using \eqref{Murata} combined with Young then H\" older inequalities
the integral above is bounded by
$$C_M\|\Omega(\cdot,s)v\|_{L^2(\mathbb{R},L^2_{-\sigma})}^2\le C_M \tilde C_K^2\|v\|_{L^2_x}.$$
where, for the last inequality we employed \eqref{est:katoo}.
Consequently, there exist a constant $C_2$ such that for any $t,\
s\in\mathbb{R}:$
$$\|T(t,s)v\|_{L^2_x}\le C_2\|v\|_{L^2_x}.$$
This finishes the proof of the Lemma. $\Box$

Fix now $2<p_0<\infty$ and let $p'_0=p_0/(p_0-1).$ By applying
Riesz-Thorin interpolations to the operators $T(t,s)$ satisfying for
all $t,\ s\in\mathbb{R}:$
 \begin{eqnarray}
 \|T(t,s)\|_{L^2\mapsto L^2}&\le &C_2\nonumber\\
 \|T(t,s)\|_{L^1\cap L^{p'_0}\mapsto
L^{p_0}}&\le &\frac{C_{p_0}\log(2+|t-s|)}
{(1+|t-s|)^{1-\frac{2}{p_0}}}\nonumber
 \end{eqnarray}
we obtain that for any $2\le p\le p_0$ there exists a constant
$C_{p_0,p}$ between $C_2$ and $C_{p_0}$ such that:
$$\|T(t,s)\|_{L^{q'}\cap L^{p'}\mapsto
L^{p}}\le \frac{C_{p_0,p}\log^{\frac{1-2/p}{1-2/p_0}}(2+|t-s|)}
{(1+|t-s|)^{1-\frac{2}{p}}},\ {\rm where}\ p'=\frac{p}{p-1},\
 q'=p'\frac{p_0-2}{p_0-p'}.$$ Finally, using \eqref{def:T} and the
estimates for the Schr\" odinger group \eqref{est:Lp} we get:

\begin{theorem}\label{th:lin2} Fix $2<p_0<\infty$ and assume that
$\|<x>^{4\sigma /3}\psi_E(t)\|_{H^2}<\varepsilon_1,\ t\in\mathbb{R}$
where $\varepsilon_1$ is the constant obtained in
Theorem~\ref{th:lin1}. Then there exists the constants $C_2,\
C_{p_0,p}>0$ such that for all $2\le p\le p_0$ and $t,\
s\in\mathbb{R}$ the following estimates hold:
\begin{eqnarray}
\|\Omega(t,s)\|_{L^2\mapsto L^2}&\le &C_{2};\nonumber\\
\|\Omega(t,s)\|_{L^{q'}\cap L^{p'}\mapsto L^p}&\le &
\frac{C_{p_0,p}\log(2+|t-s|)^{\frac{1-2/p}{1-2/p_0}}}{|t-s|^{1-\frac{2}{p}}},
\ {\rm where}\ p'=\frac{p}{p-1},\
 q'=p'\frac{p_0-2}{p_0-p'}.\nonumber
\end{eqnarray}
\end{theorem}

Note that the estimates for the family of operators $\Omega(t,s)$
given by the above theorem are similar to the standard
$L^{p'}\mapsto L^p$ estimates for Schr\" odinger operators
(\ref{est:Lp}) except for the logarithmic correction and a smaller
domain of definition $L^{q'}\cap L^{p'}\subset L^{p'}$ where $q'<p'$
when $p'<2.$ If we would have proven \eqref{est:Timp} then we could
use $p_0=\infty,$ hence $q'=p'$ in the above theorem and obtain:
$$\|\Omega(t,s)\|_{L^{p'}\mapsto L^p}\le
\frac{C_{p}\log(2+|t-s|)^{1-2/p}}{|t-s|^{1-\frac{2}{p}}} \qquad {\rm
where}\ p'=\frac{p}{p-1}.$$

\bigskip
\noindent{\bf Acknowledgements:} E. Kirr was partially supported by
NSF grants DMS-0405921, DMS-0603722 and DMS-0707800.

\bibliographystyle{plain}
\bibliography{ref}

\end{document}